\documentclass[12pt]{article}
\usepackage{graphicx}
\usepackage{amsmath,amsthm,amssymb,enumerate}
\usepackage{BOONDOX-cal}
\usepackage{url}
\usepackage{euscript,mathrsfs}
\usepackage{bm}
\usepackage[left=1.5cm,right=1.5cm,top=3.5cm,bottom=3.5cm]{geometry}
\usepackage{color}
\catcode`\@=11 \@addtoreset{equation}{section}

\catcode`\@=12
\allowdisplaybreaks

\newtheorem{Theorem}{Theorem}[section]
\newtheorem{Proposition}[Theorem]{Proposition}
\newtheorem{Lemma}[Theorem]{Lemma}
\newtheorem{Corollary}[Theorem]{Corollary}

\theoremstyle{definition}
\newtheorem{Definition}[Theorem]{Definition}

\newtheorem{Remark}[Theorem]{Remark}

\newcommand{\bTheorem}[1]{
\begin{Theorem} \label{T#1} }
\newcommand{\eT}{\end{Theorem}}

\newcommand{\bProposition}[1]{
\begin{Proposition} \label{P#1}}
\newcommand{\eP}{\end{Proposition}}

\newcommand{\bLemma}[1]{
\begin{Lemma} \label{L#1} }
\newcommand{\eL}{\end{Lemma}}

\newcommand{\bCorollary}[1]{
\begin{Corollary} \label{C#1} }
\newcommand{\eC}{\end{Corollary}}

\newcommand{\bRemark}[1]{
\begin{Remark} \label{R#1} }
\newcommand{\eR}{\end{Remark}}

\newcommand{\bDefinition}[1]{
\begin{Definition} \label{D#1} }
\newcommand{\eD}{\end{Definition}}

\newcommand{\vcg}[1]{{\pmb #1}}

\newcommand{\Del}{\Delta_x}

\newcommand{\bu}{\mathbf u}

\newcommand{\vrd}{\varrho_{\delta}}
\newcommand{\vud}{\vu_{\delta}}
\newcommand{\Zd}{Z_{\delta}}

\newcommand{\bfphi}{\boldsymbol{\varphi}}

\newcommand{\xxi}{\boldsymbol{\xi}}

\newcommand{\bFormula}[1]{
\begin{equation} \label{#1}}
\newcommand{\eF}{\end{equation}}

\newcommand{\Ov}[1]{\overline{#1}}

\newcommand{\vr}{\varrho}
\newcommand{\vre}{\vr_\ep}

\newcommand{\vue}{\vu_\ep}

\newcommand{\vu}{\vc{u}}
\newcommand{\vc}[1]{{\bf #1}}

\newcommand{\Div}{{\rm div}_x}
\newcommand{\Grad}{\nabla_x}

\newcommand{\tn}[1]{\mathbb{#1}}
\newcommand{\dx}{{\rm d} {x}}

\newcommand{\dt}{{\rm d} t }

\newcommand{\dxdt}{\dx\dt}
\newcommand{\intO}[1]{\int_{\Omega} #1 \ \dx}

\newcommand{\intTO}[1]{\int_0^T\!\!\! \int_{\Omega} #1 \ \dxdt}

\newcommand{\ep}{\varepsilon}

\newcommand{\R}{\mathbb{R}}

\newcommand{\pt}{\partial_{t}}
\newcommand{\vrs}{\vr^*}
\newcommand{\eq}[1]{\begin{equation}
\begin{split}
#1
\end{split}
\end{equation}}
\newcommand{\eqh}[1]{\begin{equation*}
\begin{split}
#1
\end{split}
\end{equation*}}
\newcommand{\lr}[1]{\left( #1 \right)}
\newcommand{\vS}{\mathbb{S}}
\newcommand{\vw}{\vc{w}}

\newcommand{\intT}[1]{\int_0^T #1 \ \dt}
\newcommand{\lap}{\Delta}
\newcommand{\A}{{\cal{A}}}

\definecolor{Cgrey}{rgb}{0.85,0.85,0.85}
\definecolor{Cblue}{rgb}{0.50,0.85,0.85}
\definecolor{Cred}{rgb}{1,0,0}
\definecolor{fancy}{rgb}{0.10,0.85,0.10}

\newcommand\Cbox[2]{%
    \newbox\contentbox%
    \newbox\bkgdbox%
    \setbox\contentbox\hbox to \hsize{%
        \vtop{
            \kern\columnsep
            \hbox to \hsize{%
                \kern\columnsep%
                \advance\hsize by -2\columnsep%
                \setlength{\textwidth}{\hsize}%
                \vbox{
                    \parskip=\baselineskip
                    \parindent=0bp
                    #2
                }%
                \kern\columnsep%
            }%
            \kern\columnsep%
        }%
    }%
    \setbox\bkgdbox\vbox{
        \color{#1}
        \hrule width  \wd\contentbox %
               height \ht\contentbox %
               depth  \dp\contentbox
        \color{black}
    }%
    \wd\bkgdbox=0bp%
    \vbox{\hbox to \hsize{\box\bkgdbox\box\contentbox}}%
    \vskip\baselineskip%
}

\date{}

\begin{document}


\title{Two-phase compressible/incompressible Navier--Stokes system with inflow-outflow boundary conditions}

\author{Milan Pokorn\'y
\thanks{The work of M. P. has been supported by the Czech Science Foundation, project No. 19-04243S}
\and
Aneta Wr\'oblewska-Kami\'nska
\thanks{The work of A. W-K. has been supported by the grant of National Science Centre Poland Sonata Bis UMO-2020/38/E/ST1/00469.}
 \and 
Ewelina Zatorska
\thanks{The research of E. Z. leading to these results has been funded by the EPSRC Early Career Fellowship no. EP/V000586/1.
This work was also partially supported by the Simons Foundation Award No 663281 granted to the Institute of Mathematics of the Polish Academy of Sciences for the years 2021-2023.
}
}

\date{\today}

\maketitle

\bigskip

{
\footnotesize
\centerline{$^*\;$Charles University, Faculty of Mathematics and Physics, Mathematical Institute}
\centerline{of Charles University, Sokolovsk\'a 83, 186 75 Praha 8, Czech Republic}
\centerline{\small \texttt{pokorny@karlin.mff.cuni.cz}}

\bigbreak
\centerline{$^\dagger\;$Institute of Mathematics, Polish Academy of Sciences}
\centerline{ul. \'Sniadeckich 8,
00-656 Warszawa, Poland}
\centerline{\small \texttt{awrob@impan.pl}}

\bigbreak
\centerline{$^\ddagger\;$Department of Mathematics, Imperial College London}
\centerline{South Kensington Campus -- SW7 2AZ, London, UK}
\centerline{\small \texttt{e.zatorska@imperial.ac.uk}}

}

\bigbreak

\begin{abstract}
We prove the existence of a weak solution to the compressible Navier--Stokes system with 
singular pressure that explodes when density achieves its congestion level. This is a quantity whose initial value evolves according to the transport equation. We then prove that the ``stiff pressure" limit gives rise to the two-phase compressible/incompressible system with congestion constraint describing the free interface. We prescribe of the velocity at the boundary and the value of density at the inflow part of the boundary of a general bounded $C^2$ domain.
There are no restrictions on the size of the boundary conditions.
\end{abstract}



\bigskip

\section{Introduction}
\label{i}
The fluid-type equations are often used as macroscopic models for collective dynamics. In the present paper we are particularly interested in a system that has been analysed in \cite{DeMiZa, DeMiNaZa} as a model the motion of big crowds. It is a two phase compressible-incompressible model describing the evolution of some averaged macroscopic quantities describing the crowd: the velocity $\vu$, the density $\vr$ and the congestion density $\vr^*$. The latter describes the preferences of the individuals, or their physical dimensions, that restrict the neighbours from being too close to each other (from penetrating each other). It is set for each individual in the crowd at the initial time, and simply transported along with the flow. This means that when the density of the crowd $\vr$ reaches this constraining value $\vr^*$, the crowd behaves like the incompressible fluid. When the density  $\vr$ is strictly less than $\vr^*$ the crowd behaves like compressible fluid, except that the particles move freely as there is no contact between them. The behaviour of the crowd is described by either incompressible or compressible Navier-Stokes equations on the moving subdomains separated by the interphase described by the relation:
\eq{\label{condition}
\pi(\vr^*-\vr)=0.}
Here $\pi$ is the unknown ``pressure'', or rather, the Lagrangian multiplier associated with incompressibility condition satisfied by the velocity. Relation \eqref{condition} states that $\pi$ appears only on the subdomain with congestion. For $\vr<\vr^*$, on the other hand, $\pi$ is equal to $0$.

 The similar free boundary problem was already analysed by Lions and Masmoudi in \cite{LM99} for $\vr^*=1$. The authors showed that the two-phase system can be approximated by purely compressible Navier-Stokes equations with the pressure $\pi_n(\vr)\approx\vr^{\gamma_n}$ and $\gamma_n\to\infty$. The same kind of limit passage was also investigated later on for the PDE models of tumor growth \cite{PeVa, VauZat}. In the current paper we will focus on another approximation of the unknown pressure
 $\pi_\ep\approx\ep \frac1{\vr^*-\vr}$,  which has some benefits from the numerical perspective, see \cite{DeMiNaZa,DeHuNa}. Similar forms and asymptotic limits of the singular pressure appear in the models of traffic models \cite{BeDeDeRa, BeBr, BeDeBlMoRaRo}, collective dynamics \cite{DeHuNa}, or granular flow \cite{Maury, Perrin}.

 All of the previous analytical results for the derivation of the two-phase compressible-incompressible system were obtained either for the whole space case, or for the bounded domain with zero Dirichlet boundary condition. In this paper, we want to extend the analysis to the setting where the inflow and outflow of the crowd is allowed, making the model suitable to describe various evacuation scenarios. Some numerical simulations for the hyperbolic version of such model were already performed in \cite{DeMiNaZa}.

 Our starting point is the following system of equations in $(0,T)\times \Omega$:
 \begin{subequations}\label{main}
\eq{
\label{i1}
\partial_t\vr+\Div (\vr \vu)= 0,
}
\eq{
\label{i2}
\partial_t(\vr\vu)+\Div (\vr \vu \otimes \vu) + \Grad \pi_\ep\Big(\frac{\vr}{\vrs}\Big)-\Div \mathbb{S}(\Grad \vu)=\vr(\vc{w}-\vu), }
\eq{\label{i3}
\pt \vrs+\vu\cdot\Grad\vrs=0,}
\end{subequations}
where $\Omega \subset R^d$, $d=2,3$  is a bounded domain of class $C^2$,  $\vr = \vr(t,x)$ is the unknown mass density,  $\vu = \vu(t,x)$  is
the unknown velocity, $\vrs=\vrs(t,x)$ is the unknown congestion density, $(t,x)\in (0,T)\times\Omega\equiv Q_T$, $\ep>0$, and $\vc{w}=\vw(x)$ is given.

The stress tensor $\mathbb{S}$ will be specified below (as stress tensor for compressible Newtonian fluid).
The pressure $\pi$ depends on the ratio of densities $\frac{\vr}{\vrs}$. 
It models a constraint on the density of the fluid $\vr$ that cannot exceed the value $\vrs$. More precisely, we take
\begin{equation} \label{i5a}
\pi_\ep\Big(\frac{\vr}{\vrs}\Big) =\ep\frac{\big(\frac{\vr}{\vrs}\big)^\alpha}{\big(1-\frac{\vr}{\vrs}\big)^\beta} 
\end{equation} 
with some $\alpha >1$ and $\beta>\frac 52$ id $d-3$ and $\beta>2$ if $d=3$. Note that the pressure fulfills 
\begin{equation} 
\label{i6}
\pi_\ep(0)=0, \quad \pi_\ep'(z)>0 \quad \mbox{ for } 0<z<1.
\end{equation}
This is similar to the so-called hard sphere pressure considered, e.g., in \cite{CHNY2}, that can be viewed as a special case of system  \eqref{main} with $\vrs=const$. Similarly as in \cite{CHNY2} we can relax the assumption on $\pi$ for $\vr<\vrs$ by addition of non-monotone pressure that vanishes for $\vr\to\vrs$. To avoid unnecessary technicalities, we skip this point in this paper. 
We consider our system together  with initial conditions
\begin{equation}\label{initc}
\vr(0)=\vr_0,\quad\vr\vu(0)=\vr_0\vu_0,\quad \vrs(0)=\vrs_0,
\end{equation}
and boundary conditions
\begin{equation}
\label{i4}
\vc{u} |_{\partial \Omega} = \vu_B, \quad \vr|_{\Gamma_{\rm in}} = \vr_B, \quad \vrs|_{\Gamma_{\rm in}} = \vrs_B,
\end{equation}
where
\begin{equation}
\label{i5}
\Gamma_{\rm in} = \left\{ x \in \partial \Omega \ \Big| \ {\vc{u}_B} \cdot \vc{n} < 0 \right\}, \quad 
\Gamma_{\rm out} = \left\{ x \in \partial \Omega \ \Big| \ {\vc{u}_B} \cdot \vc{n} \geq 0 \right\}
\end{equation}
(we include to the outflow part of the boundary also the part where the normal velocity component is zero). 

Note that system \eqref{i1}--\eqref{i3} reminds system studied (for homogeneous boundary conditions for the velocity) in \cite{M3NPZ}; there, the role of the congestion density was played by the entropy. Using a similar idea as in the above mentioned paper we introduce a new quantity $Z:=\frac{\vr}{\vrs}$; then the new unknown function $Z$ satisfies (at least formally, for smooth solutions) the continuity equation and we obtain the following system
 \begin{subequations}\label{main_transformed}
\eq{
\label{ia1}
\partial_t\vr+\Div (\vr \vu)= 0,
}
\eq{
\label{ia2}
\partial_t(\vr\vu)+\Div (\vr \vu \otimes \vu) + \Grad \pi_\ep(Z)-\Div \mathbb{S}(\Grad \vu)=\vr(\vc{w}-\vu), }
\eq{\label{ia3}
\pt Z+\Div (Z \vu)=0}
\end{subequations}
with initial
\begin{equation}\label{initca}
\vr(0)=\vr_0,\quad\vr\vu(0)=\vr_0\vu_0,\quad Z(0)=Z_0:=\frac{\vr_0}{\vrs_0},
\end{equation}
and boundary conditions
\begin{equation}
\label{ia4}
\vc{u} |_{\partial \Omega} = \vu_B, \quad \vr|_{\Gamma_{\rm in}} = \vr_B, \quad Z|_{\Gamma_{\rm in}} = Z_B:=\frac{\vr_B}{\vrs_B}.
\end{equation}
Note, however, that by standard techniques we can get certain "better" information on the "density" only from the pressure term, therefore, as in \cite{M3NPZ}, we will consider a certain interplay of the initial and boundary conditions for $\vr$ and $\vrs$ which leads to the fact that the boundary and initial conditions for $\vr$ are controlled by the initial and boundary conditions for $Z$ (see \eqref{data_Z}). Furthermore, we also have that the initial and boundary conditions for $Z$ belong to the interval $(0,1)$. 

Combining the approach from \cite{M3NPZ} with \cite{CHNY2} we will be able to show that under certain additional technical assumptions problem \eqref{ia1}--\eqref{ia3} with initial \eqref{initca} and boundary conditions \eqref{ia4} possesses a weak solution defined below.  We even slightly improve the result from \cite{CHNY2} in the sense that we may include for global-in-time existence result the case when the velocity flux is zero, see Remark \ref{RB1}. Next, it is possible to show that also \eqref{i1}--\eqref{i3} with initial \eqref{initc} and boundary conditions \eqref{i4} has a solution: the approach is based on suitable renormalization which allows us to return back to the unknown function $\vrs$. On the other hand, we are more interested in the limit passage $\varepsilon \to 0+$; we perform the limit passage in the formulation with the unknown function $Z$ and later return back to the formulation with the function $\vrs$. Here, we follow ideas from \cite{DeMiZa} or \cite{PeZa}. We will show that with
 $\ep\to 0^+$ the weak solutions to system \eqref{main} converge in some sense to the weak solution of the target system
 \begin{subequations}\label{target}
 \begin{equation}\label{rho}
 \pt\vr+ \Div (\vr\vu) = 0, 
 \end{equation}
\begin{equation}\label{mom}
 \partial_t (\vr\vu) + \Div (\vr\vu \otimes \vu) + \Grad\pi -\Div\vS(\Grad\vu) = \vr(\vw-\vu),
 \end{equation}
\begin{equation}\label{rho_star}
 \pt\vr^*+\vu\cdot\Grad\vr^*=0,
 \end{equation}
 \begin{equation}\label{cons0}
 0\leq \vr\leq\vrs,
 \end{equation}
\begin{equation}\label{div0}
  \Div\vu =0 \ \text{in} \ \{\vr=\vrs\},
 \end{equation}
\begin{equation}
  \pi\geq 0\ \mbox{in} \ \{\vr=\vrs\},\quad \pi= 0\ \mbox{in} \ \{\vr<\vrs\}.\label{pineq0}
  \end{equation}
 \end{subequations}
We call this system a free boundary two-phase compressible/incompressible system. To justify this name note that when $\vr=\vrs$, i.e., when the density achieves its maximal value, due to condition
\eqref{div0} the system behaves like inhomogeneous incompressible Navier-Stokes equations. When on the other hand $\vr<\vrs$, the system behaves like compressible pressureless Navier-Stokes equations with time-space variable upper bound for the density. This is one of the novelties here, as before we always had to consider the background pressure in the barotropic form.

One of the inspirations for this work was the numerical paper \cite{DeMiNaZa}, in which the inviscid variant of system \eqref{main} was considered with "do nothing" boundary condition for the velocity at the outflow part of the boundary. Even though this condition is often used in numerics, it is not very suitable for analysis, in particular for global-in-time type solutions for large data.

Even though we combine approaches from several other papers, the result itself is new and requires nontrivial extensions of previously known techniques and results.

\section{Main result}
\label{M}

In what follows, we formulate main results of the paper. Before doing so, we state the main assumptions on the data of our problem.
Concerning the given field $\vc{w}$, we consider (for simplicity)
\begin{equation} \label{As1}
\vc{w} \in L^\infty(Q_T;R^d).
\end{equation}
Further, the stress tensor $\mathbb{S}$ is characteristic for the Newtonian fluid and it is given by
\begin{equation}\label{is}
\mathbb{S}(\Grad \vu) = \mu \left( \Grad \vu + \Grad^t \vu \right) + \lambda \Div \vu \mathbb{I}, \ \mu > 0 , \ \lambda \geq 0.
\end{equation}
The pressure $\pi_\ep(\cdot)$ has the form \eqref{i5a} with $\alpha >1$ and $\beta >\frac 52$.

Similarly as in \cite{CHNY2}, we consider the following regularity assumptions
\begin{equation} \label{Ass1}
{ \vr_B, \vrs_B \in C(\Gamma_{\rm in})}, \quad \vu_B \in C^2(\partial \Omega; R^d), \quad \int_{\partial \Omega} \vu_B\cdot\vc{n}\ {\rm d}S_x \dt \geq 0.
\end{equation}

Furthermore, we assume that
\begin{equation} \label{Ass2}
\begin{array}{c}
\displaystyle 0<  \vr_0<\vrs_0,  \  \mbox{a.e. in }\Omega, \quad \vrs_0\in L^\infty(\Omega),\\
\displaystyle \intO{H_\ep\lr{\frac{\vr_0}{\vrs_0}}}<\infty,\quad {\rm ess\ inf}_\Omega\, \vr_0 >0, \quad {\rm ess\ inf}_\Omega\, (\vrs_0 -\vr_0) >0,\\
\vu_0 \in L^2(\Omega;R^d),
\end{array}
\end{equation}
where \bFormula{Hf}
\begin{split}
H_\ep(z) &= z\int_0^z\frac{\pi_\ep(s)}{s^2}\ {\rm d}s. 
\end{split}
\eF
For the boundary data,
\begin{equation} \label{Ass3}
\displaystyle {0<  \vr_B <\vrs_B,  \  \mbox{a.e. on }\Gamma_{\rm in}}, \quad  {\rm ess\ inf}_{\Gamma_{\rm in}}\, \vr_B >0, \quad {\rm ess\ inf}_{\Gamma_{\rm in}}\, (\vrs_B -\vr_B) >0.
\end{equation}
Note that these assumptions yield that the initial energy
$$
E_0 := \intO{\Big(\frac 12 \vr_0 |\vu_0|^2 + H_\ep\Big(\frac{\vr_0}{\vrs_0}\Big)\Big)} <+\infty
$$
as well as that there are positive constants $c_*$ and $c^*$ such that
$$
c_* \leq \frac{1}{\vrs_0} \leq c^*, \quad c_*\leq \frac{1}{\vrs_B} \leq c^* \quad \text{a.e.}
$$
Whence, rewriting our problem to the form \eqref{ia1}--\eqref{ia4}, we immediately have that
\begin{equation} \label{data_Z}
{ 
c_* \vr_0 \leq Z_0 \leq c^* \vr_0 \mbox{ a.e. in }\Omega,
\quad c_* \vr_B\leq Z_B \leq c^*\vr_B \quad \text{ a.e. on }\Gamma_{\rm in}
}
\end{equation}

Let us now introduce the definition of a weak solution to problem \eqref{ia1}--\eqref{ia4}.
\begin{Definition} \label{def_in_Z}
We say that $(\vr, \vu, Z)$ is a bounded energy weak solution of problem \eqref{ia1}--\eqref{ia4}  on a time interval $(0,T)$ if  the following five conditions are satisfied.

\begin{description}
\item{1.} 
The triple of functions $(\vr, \vu, Z)$ fulfills:
\eq{\label{fs-} 
&0 \leq c_\star \vr \leq Z \leq c^\star \vr \mbox{ a.e. in } Q_T,  \quad \mbox{for}\quad 0<c_\star\leq c^\star <\infty,\\
 &0 \leq Z<1  \text{ a.e. in } (0,T) \times \Omega, \quad  \pi_\ep(Z)\in L^1(0,T; L^1(\Omega))
\\
&\vu \in L^2(0,T; W^{1,2}(\Omega; R^d)), \quad 
\vu|_{I\times\partial \Omega} = \vu_B.
}
\item{2.}
{ The function $\vr\in C_{\rm weak}([0,T], L^1(\Omega))$
satisfies the integral identity
\bFormula{ce-}
\begin{split}
&\intO{\vr(\tau,\cdot)\varphi(\tau,\cdot)} - \intO{\vr_0(\cdot)\varphi(0,\cdot)} \\
&= \int_0^\tau\intO{\Big(\vr\partial_t\varphi + \vr \vu \cdot \Grad \varphi\Big) }{\rm d}t - \int_0^\tau\int_{\Gamma_{\rm in}} \vr_B \vu_B \cdot \vc{n} \varphi 
\ {\rm d}S_x{\rm d} t 
\end{split}
\eF
for any $\tau\in [0,T]$ and $\varphi \in C_c^1([0,T]\times({\Omega}\cup\Gamma_{\rm in}))$.}
\item{3.} 
{ The function $\vr\vu\in  C_{\rm weak}([0,T], L^{1}(\Omega;R^d))$ satisfies the integral identity
\bFormula{me-}
\begin{split}
&\intO{\vr\vu(\tau,\cdot)\cdot\bfphi(\tau,\cdot)} - \intO{\vr_0\vu_0(\cdot)\cdot\bfphi(0,\cdot)} \\
&= \int_0^\tau \intO{\Big( \vr\vu\cdot\partial_t\bfphi+\vr \vu \otimes \vu : \Grad \bfphi + \pi_\ep(Z)\Div \bfphi - \mathbb{S}(\Grad \vu) : \Grad \bfphi \Big)}{\rm d}t\\
&\quad+\int_0^\tau \intO{\vr(\vw-\vu)\cdot\bfphi}{\rm d}t
\end{split}	
\eF
for any  $\tau\in [0,T]$ and any $\bfphi \in C^1_c ([0,T]\times\Omega; R^d)$.}
\item{4.}
{ The function $Z\in C_{\rm weak}([0,T], L^1(\Omega))$
satisfies the integral identity 
\bFormula{ceZ-}
\begin{split}
&\intO{Z(\tau,\cdot)\varphi(\tau,\cdot)} - \intO{Z_0(\cdot)\varphi(0,\cdot)} \\
&= \int_0^\tau\intO{\Big(Z\partial_t\varphi + Z \vu \cdot \Grad \varphi\Big) }{\rm d}t - \int_0^\tau\int_{\Gamma_{\rm in}} Z_B \vu_B \cdot \vc{n} \varphi 
\ {\rm d}S_x{\rm d} t 
\end{split}
\eF
for any $\tau\in [0,T]$ and $\varphi \in C_c^1([0,T]\times({\Omega}\cup\Gamma_{\rm in}))$.}
\item{5.} 
{ There is a Lipschitz extension $\vu_\infty\in W^{1,\infty}(\Omega;R^d)$ of the vector field $\vu_B$ such that the following energy inequality holds 
\bFormula{ei-}
\begin{split}
&\intO{\Big(\frac 12\vr|\vu-\vu_\infty|^2+H_\ep(Z)\Big)(\tau)}+\int_0^\tau\intO{\tn S(\Grad(\vu-\vu_\infty)):\Grad(\vu-\vu_\infty)}{\rm d}t \\
&\quad 
+\int_0^\tau\int_\Omega\pi_\ep(Z) \Div\vu_\infty\ {\rm d}x{\rm d}t \\
&\le \intO{\Big(\frac 12\vr_0|\vu_0-\vu_\infty|^2+H_\ep(Z_0)\Big)}
-\int_0^\tau\intO{\vr\vu\cdot\Grad\vu_\infty\cdot(\vu-\vu_\infty)}{\rm d}t \\
&\quad -\int_0^\tau\intO{\tn S(\Grad\vu_\infty):\Grad(\vu-\vu_\infty)}{\rm d}t- \int_0^\tau\int_{\Gamma_{\rm in}} H_\ep(Z_B)\vu_B\cdot\vc n \ {\rm d}S_x{\rm d}t \\
&\quad + \int_0^\tau\intO{\vr(\vc{w}-\vu)\cdot(\vu-\vu_\infty)}{\rm d}t
\end{split}	
\eF
for a.a. $\tau\in (0,T)$. }
\end{description}
\end{Definition}
\noindent
{
\bRemark{R1-}
The continuity equations \eqref{ia1} and \eqref{ia3} give rise to
\bFormula{mi}
\intO{X(\tau)}\le\intO{X_0}-\int_0^\tau\int_{\Gamma_{\rm in}}X_B\vu_B\cdot\vc n\ {\rm d} S_x{\rm d}t
\eF
for all $X=\vr,Z$, and  $\tau\in [0,T]$. It can be obtained by taking in \eqref{ce-} and \eqref{ceZ-}  test functions $\varphi=\varphi_\eta$, where 
$$\varphi_\eta(x) = \left\{\begin{array}{rl} 1 & \text{if } {\rm dist}\,(x,\Gamma_{\rm out})>\eta \\
\frac{1}{\eta} {\rm dist}\,(x,\Gamma_{\rm out}) & \text{if } {\rm dist}\,(x,\Gamma_{\rm out}) \leq \eta,
\end{array}
\right.
$$ 
and then by letting $\eta\to 0+$.
\eR

\begin{Definition} \label{def_renor}
We call $(Z,\vu)$ a renormalized solution of the continuity equation \eqref{ia3} provided:
\begin{itemize}
\item {$Z\in L^\infty(0,T;L^\infty(\Omega))$}, and $\vu\in  L^2(0,T; W^{1,2}(\Omega,R^3))$,
\item $(Z,\vu)$ satisfies  the weak formulation of the continuity equation \eqref{ceZ-},
\item for any $b\in C^1[0,1]$, $b(Z)\in C_{\rm weak}([0,T];L^1(\Omega))$
the weak  formulation of the renormalized equation is satisfied, i.e.,
\begin{equation} 
\label{P3}
\begin{split}
&\intO{ b(Z)  (\tau,\cdot)\varphi(\tau,\cdot)} -\intO{ b(Z_0) (\cdot)\varphi(0,\cdot)} \\
&= 
{
\int_0^\tau\intO{ \Big({{ b(Z)\partial_t\varphi+} b(Z) \vu \cdot \Grad \varphi 
-\varphi\left( b'(Z) Z - b(Z) \right) \Div \vu \Big)}} {\rm d}t
}\\  
&\quad - \int_0^\tau\int_{\Gamma_{\rm in}} b(Z_B) 
\vu_B \cdot \vc{n} \varphi \ {\rm d}S_x {\rm d} t 
\end{split}
\end{equation}
for any $\varphi \in C^1_c([0,T]\times({\Omega}\cup\Gamma_{\rm in}))$.
\end{itemize}
A weak solution to problem \eqref{ia1}--\eqref{ia4} satisfying in addition  renormalized continuity equations \eqref{P3} for both $(Z,\bu)$ and $(\vr,\vu)$ is called a {\it renormalized weak solution}.
\end{Definition}

\bRemark{R1aa} 
Note that due to \cite[Lemma 3.1]{CHJN}, and since we know that since both $\vr$ and $Z$ are essentially bounded functions and $\vu \in L^2(0,T;W^{1,2}(\Omega;R^d))$, any weak solution in the sense of Definition \ref{def_in_Z} is in fact a renormalized weak solution, provided $\Gamma_{\rm in}$ is a $C^2$ open $(d-1)$ dimensional manifold. 
\eR

\bigskip

In the proof that follows we will need that the boundary data for the velocity, $\vu_B$, can be extended to the whole $\Omega$ in such a way that the extension is sufficiently smooth and its divergence is non-negative. This follows from the following result (see \cite[Lemmata 5.1, 5.2, and 5.3]{CHNY2}

\begin{Lemma} \label{Lextrem}
 Let $\vc V \in  C^2(\partial\Omega;R^d)$ 
be a  vector field on the boundary  $\partial\Omega$ of a bounded $C^2$ domain $\Omega$. Let 
$$
\int_{\partial \Omega} \vc V \cdot \vc{n} \ {\rm d}S_x \geq 0.
$$
Then there exist 
a vector field
{ \bFormula{Le}
\vc V_\infty\in W^{2,q}(\Omega;R^d), \quad 1\leq q<\infty, \quad \Div {\vc V}_\infty\ge 0 \quad \text{ a.e. in } \Omega
\eF}
verifying ${\vc V}_\infty|_{\partial\Omega}=\vc V$. \\
If in addition
$$
\int_{\partial \Omega} \vc V \cdot \vc{n} \ {\rm d}S_x =K> 0,
$$
then the extension $\vc V_\infty$ satisfy
{ \bFormula{Lep}
 \Div {\vc V}_\infty\ge 0 \quad \text{ a.e. in } \Omega,\qquad and\qquad {\rm ess}\inf_{{\cal O}}\lr{\Div {\vc V}_\infty}\ge C> 0,
\eF}
where ${\cal O}$ is an open subset satisfying $\bar{\cal O}\subset\Omega$.
\end{Lemma} 

Our first result is  a global-in-time existence theorem for solutions defined above.
\begin{Theorem} 
\label{TM1!}
Let $\Omega \subset R^d$, $d = 2,3$, be a bounded domain of class $C^{2}$ such that $\Gamma_{\rm in}$ is an open $C^2$ $d-1$ dimensional manifold. Let $\ep>0$, $T>0$.
{
Under the assumptions \eqref{As1}--\eqref{Ass3} the problem \eqref{main_transformed}--\eqref{ia4}
admits at least one bounded energy weak solution $(\vr, \vu, Z)$ on $(0,T)$ in the sense of Definition~\ref{def_in_Z}. Moreover  $(\vr,\vu,Z)$ is a renormalized solution in the sense of Definition \ref{def_renor}.
}
\end{Theorem}

Next we consider solutions to problem {
\eqref{main}--\eqref{i5} with \eqref{is}}. 

\begin{Definition} \label{def_in_vrs}
We say that $(\vr, \vu, \vrs)$ is a bounded energy weak solution of problem {
\eqref{main}--\eqref{i5} with \eqref{is}}  on a time interval $(0,T)$ if  the following five conditions are satisfied.

\begin{description}
\item{1.} 
The triple of functions belongs fulfills:
\eq{\label{gs-} 
&0 \leq  \vr < \vr^* \mbox{ a.e. in } Q_T\\
 & c_* \leq \frac{1}{\vrs} \leq c^*  \text{ a.e. in } (0,T) \times \Omega, \quad  \pi_\ep\Big(\frac{\vr}{\vrs}\Big)\in L^1(0,T; L^1(\Omega))
\\
&\vu \in L^2(0,T; W^{1,2}(\Omega; R^d)), \quad 
\vu|_{I\times\partial \Omega} = \vu_B.
}
\item{2.}
{ The function $\vr\in C_{\rm weak}([0,T], L^1(\Omega))$
satisfies the integral identity
\bFormula{ge-}
\begin{split}
&\intO{\vr(\tau,\cdot)\varphi(\tau,\cdot)} - \intO{\vr_0(\cdot)\varphi(0,\cdot)} \\
&= \int_0^\tau\intO{\Big(\vr\partial_t\varphi + \vr \vu \cdot \Grad \varphi\Big) }{\rm d}t - \int_0^\tau\int_{\Gamma_{\rm in}} \vr_B \vu_B \cdot \vc{n} \varphi 
\ {\rm d}S_x{\rm d} t 
\end{split}
\eF
for any $\tau\in [0,T]$ and $\varphi \in C_c^1([0,T]\times({\Omega}\cup\Gamma_{\rm in}))$.}
\item{3.} 
{ The function $\vr\vu\in  C_{\rm weak}([0,T], L^{1}(\Omega;R^d))$ satisfies the integral identity
\bFormula{gge-}
\begin{split}
&\intO{\vr\vu(\tau,\cdot)\cdot\bfphi(\tau,\cdot)} - \intO{\vr_0\vu_0(\cdot)\bfphi(0,\cdot)} \\
&= \int_0^\tau \intO{\Big( \vr\vu\cdot\partial_t\bfphi+\vr \vu \otimes \vu : \Grad \bfphi + \pi_\ep \Big(\frac{\vr}{\vrs}\Big)\Div \bfphi - \mathbb{S}(\Grad \vu) : \Grad \bfphi \Big)}{\rm d}t\\
&\quad+\int_0^\tau \intO{\vr(\vw-\vu)\cdot\bfphi}{\rm d}t
\end{split}	
\eF
for any  $\tau\in [0,T]$ and any $\bfphi \in C^1_c ([0,T]\times\Omega; R^d)$.}
\item{4.}
{ The function $\vrs\in C_{\rm weak}([0,T], L^1(\Omega))$
satisfies the integral identity 
\bFormula{geZ-}
\begin{split}
&\intO{\vrs(\tau,\cdot)\varphi(\tau,\cdot)} - \intO{\vrs_0(\cdot)\varphi(0,\cdot)} \\
&= \int_0^\tau\intO{\Big(\vrs\partial_t\varphi + \vrs \Div (\varphi \vu) \Big) }{\rm d}t - \int_0^\tau\int_{\Gamma_{\rm in}} \vrs_B \vu_B \cdot \vc{n} \varphi 
\ {\rm d}S_x{\rm d} t 
\end{split}
\eF
for any $\tau\in [0,T]$ and $\varphi \in C_c^1([0,T]\times({\Omega}\cup\Gamma_{\rm in}))$.}
\item{5.} 
{ There is a Lipschitz extension $\vu_\infty\in W^{1,\infty}(\Omega;R^d)$ of the vector field $\vu_B$ such that the following energy inequality holds 
\bFormula{gi-}
\begin{split}
&\intO{\Big(\frac 12\vr|\vu-\vu_\infty|^2+H_\ep\Big(\frac{\vr}{\vrs}\Big)\Big)(\tau)}+\int_0^\tau\intO{\tn S(\Grad(\vu-\vu_\infty)):\Grad(\vu-\vu_\infty)}{\rm d}t \\
&\quad 
+\int_0^\tau\int_\Omega\pi_\ep\Big(\frac{\vr}{\vrs}\Big)\Div\vu_\infty \ {\rm d}x{\rm d}t \\
&\le \intO{\Big(\frac 12\vr_0|\vu_0-\vu_\infty|^2+H_\ep\Big({
\frac{\vr_0}{\vr^*_0}}
\Big)\Big)}
-\int_0^\tau\intO{\vr\vu\cdot\Grad\vu_\infty\cdot(\vu-\vu_\infty)}{\rm d}t \\
&\quad -\int_0^\tau\intO{\tn S(\Grad\vu_\infty):\Grad(\vu-\vu_\infty)}{\rm d}t- \int_0^\tau\int_{\Gamma_{\rm in}} H_\ep\Big(\frac{\vr_B}{\vrs_B}\Big)\vu_B\cdot\vc n \ {\rm d}S_x{\rm d}t \\
&\quad + \int_0^\tau \intO{\vr (\vc{w}-\vu)\cdot (\vu-\vu_\infty)}{\rm d}t
\end{split}	
\eF
for a.a. $\tau\in (0,T)$. }
\end{description}

{
Bounded energy weak solution to problem \eqref{main}--\eqref{i5} with \eqref{is} satisfying in addition  renormalized continuity equation \eqref{P3} for $(\vr,\vu)$ is called a  renormalized weak solution.
}
\end{Definition}


We have the following result.

\begin{Theorem}\label{TM2}
Let $\Omega \subset R^d$, $d = 2,3$, be a bounded domain of class $C^{2}$ such that $\Gamma_{\rm in}$ is an open $C^2$ $d-1$ dimensional manifold. Let $\ep>0$, $T>0$.
{
Under the assumptions \eqref{As1}--\eqref{Ass3} the problem \eqref{main}--\eqref{i5} with \eqref{is}
admits at least one renormalized bounded energy weak solution $(\vr, \vu, \vrs )$ on $(0,T)$ in the sense of Definition~\ref{def_in_vrs}.
}
\end{Theorem}

Our next results concern the limit passage $\ep\to 0$. We will show that when $\ep\to0$ the weak solutions from the previous section approximate weak solutions to the system \eqref{target} defined below.
\begin{Definition}\label{Def:limit}
A quadruple $(\vr,\vu,\vrs,\pi)$ is called a global finite energy weak solution  to \eqref{target}, \eqref{is}, with the initial data \eqref{initc}, \eqref{Ass2}, and the boundary conditions \eqref{i4}, \eqref{Ass1} if for any $T>0$:
\begin{itemize}
\item[(i)] There holds:
$$0\leq\vr\leq\vrs\quad a.e.\ in \ (0,T)\times\Omega,$$ 
\eq{\label{cond:u}
\Div\vu=0\quad a.e. \ in\  \{\vr=\vrs\},} 
\eq{\label{cond:pi}
(\vr^*-\vr)\pi=0, }
and
\begin{align*}
&\vr\in C_w([0,T];L^\infty(\Omega)), \\
&\vrs\in C_w([0,T];L^\infty(\Omega)),\\
& \vu \in L^2(0,T;W^{1,2}(\Omega, R^d)),\quad \vr|\vu|^2 \in L^{\infty}(0,T; L^1(\Omega)),\\
&\pi\in {\cal M}^+ ((0,T)\times \Omega).
\end{align*}
\item[(ii)] For any $0\leq \tau\leq T$, equations \eqref{rho}, \eqref{mom}, \eqref{rho_star} are satisfied in the weak sense, more precisely:\\
-the continuity equation:
\bFormula{ge-lim}
\begin{split}
&\intO{\vr(\tau,\cdot)\varphi(\tau,\cdot)} - \intO{\vr_0(\cdot)\varphi(0,\cdot)} \\
&= \int_0^\tau\intO{\Big(\vr\partial_t\varphi + \vr \vu \cdot \Grad \varphi\Big) }{\rm d}t - \int_0^\tau\int_{\Gamma_{\rm in}} \vr_B \vu_B \cdot \vc{n} \varphi 
\ {\rm d}S_x{\rm d} t 
\end{split}
\eF
holds for any $\tau\in [0,T]$ and $\varphi \in C_c^1([0,T]\times({\Omega}\cup\Gamma_{\rm in}))$,\\
-the momentum equation:
\bFormula{gge-lim}
\begin{split}
&\intO{\vr\vu(\tau,\cdot)\cdot\bfphi(\tau,\cdot)} - \intO{\vr_0\vu_0(\cdot)\bfphi(0,\cdot)} \\
&= \int_0^\tau \intO{\Big( \vr\vu\cdot\partial_t\bfphi+\vr \vu \otimes \vu : \Grad \bfphi + \pi\Div \bfphi - \mathbb{S}(\Grad \vu) : \Grad \bfphi \Big)}{\rm d}t\\
&\quad+\int_0^\tau \intO{\vr(\vw-\vu)\cdot\bfphi}{\rm d}t
\end{split}	
\eF
holds for any  $\tau\in [0,T]$ and any $\bfphi \in C^1_c ([0,T]\times\Omega; R^d)$,\\
- the transport equation for $\vr^*$:
\bFormula{geZ-lim}
\begin{split}
&\intO{\vrs(\tau,\cdot)\varphi(\tau,\cdot)} - \intO{\vrs_0(\cdot)\varphi(0,\cdot)} \\
&= \int_0^\tau\intO{\Big(\vrs\partial_t\varphi + \vrs \Div (\varphi \vu) \Big) }{\rm d}t - \int_0^\tau\int_{\Gamma_{\rm in}} \vrs_B \vu_B \cdot \vc{n} \varphi 
\ {\rm d}S_x{\rm d} t 
\end{split}
\eF
for any $\tau\in [0,T]$ and $\varphi \in C_c^1([0,T]\times({\Omega}\cup\Gamma_{\rm in}))$.
\item[(iii)] There is a Lipschitz extension $\vu_\infty\in W^{1,\infty}(\Omega;R^d)$ of the vector field $\vu_B$ such that the following energy inequality holds 
\bFormula{gi-lim}
\begin{split}
&\intO{\frac 12\vr|\vu-\vu_\infty|^2(\tau)}+\int_0^\tau\!\!\intO{\tn S(\Grad(\vu-\vu_\infty)):\Grad(\vu-\vu_\infty)}{\rm d}t 
+\int_0^\tau\!\!\int_\Omega\pi\Div\vu_\infty \ {\rm d}x{\rm d}t \\
&\le \intO{\frac 12\vr_0|\vu_0-\vu_\infty|^2}
-\int_0^\tau\!\!\intO{\vr\vu\cdot\Grad\vu_\infty\cdot(\vu-\vu_\infty)}{\rm d}t \\
&\quad -\int_0^\tau\!\!\intO{\tn S(\Grad\vu_\infty):\Grad(\vu-\vu_\infty)}{\rm d}t + \int_0^\tau\!\! \intO{\vr (\vc{w}-\vu)\cdot (\vu-\vu_\infty)}{\rm d}t
\end{split}	
\eF
for a.a. $\tau\in (0,T)$. 
\end{itemize}
\end{Definition}
\begin{Remark}
In the above definition all the terms must make sense, in particular, $\pi$ is not only a measure, but it is sufficiently regular so that the condition \eqref{cond:pi} makes sense.
\end{Remark}

Our main theorem in this parts reads as follows.
\begin{Theorem} 
\label{TM3}
Let $\Omega \subset R^d$, $d = 2,3$, be a bounded domain of class $C^{2}$ such that $\Gamma_{\rm in}$ is an open $C^2$ $d-1$ dimensional manifold. Let $T>0$, and let assumptions \eqref{As1}--\eqref{Ass3} be satisfied.\\
If, in addition either
\eq{
\int_{\partial\Omega}\vu_B\cdot\vc{n}\, {\rm d} S_x=K>0,
}
or
\eq{
\intO{Z_0}+T\int_{\Gamma_{in}}Z_B|\vu_B\cdot\vc{n}|\, {\rm d} S_x<|\Omega|,
}
then the problem \eqref{target}, with \eqref{is}
admits at least one renormalized bounded energy weak solution $(\vr, \vu, \vrs,\pi )$ on $(0,T)$ in the sense of Definition~\ref{Def:limit}.

\end{Theorem}

The paper is organized as follows. In Section \ref{Sec:3}  we present the approximate scheme starting from the level of truncations of singular pressure at the level described by parameter $\delta$. Later on, in Section \ref{4} we obtain the uniform estimates with respect to $\delta$ and pass to the limit. As an outcome of this section we prove the Theorem \ref{TM1!} and also Theorem \ref{TM2}. In Section \ref{Sec:lim}} we recall uniform estimates with respect to $\ep$ and perform the limit passage $\ep\to0$ and conclude the proof of Theorem \ref{TM3}.

\section{Approximate solution} \label{Sec:3}
The purpose of this section is to construct approximate solutions to system \eqref{ia1}--\eqref{ia4}. We are not going to explain the details of the whole procedure but only to summarize it and to explain how can the existing literature be employed.
We approximate the singular pressure in system \eqref{main_transformed} by the truncation
\[
\pi_\delta (Z)= \begin{cases}
\pi_\ep (Z) & \text{ if } Z \in [0,1-\delta]  \\
\pi_\ep(1-\delta) + \ep(Z-1+\delta)_+^\gamma  
& \text{ if } Z \in (1- \delta,\infty),
\end{cases}
\]
where the exponent $\gamma$ has to be chosen sufficiently large in order to obtain sufficient estimates, in particular we need $\gamma>d$. 
This truncation allows to combine the arguments from \cite{M3NPZ, CHJN} in order to construct the solutions to system \eqref{main_transformed} with $\pi_\ep$ replaced by $\pi_\delta$. This consists of regularising both equations for $\vr$ and $Z$ by adding small viscosity term, we consider
\begin{equation}
\label{A1}
\partial_t\vr-\eta \Del \vr  + \Div (\vr \vu) = 0,
\end{equation}
\begin{equation}
\label{A2}
\vr(0,x)=\vr_0(x),\;\left( -\eta \ \Grad \vr + \vr \vu \right) \cdot \vc{n}|_{I\times\partial \Omega} = \left\{
\begin{array}{l} \vr_B \vu_B \cdot \vc{n} \ \mbox{if}\ [\vu_B \cdot \vc{n}](x) \leq 0,\,x\in \partial\Omega,\\
\vr \vu_B \cdot \vc{n}\;\mbox{if}\ [\vu_B \cdot \vc{n}](x) > 0,\,x\in \partial\Omega,\end{array} \right.
\end{equation}
\begin{equation}
\label{A3}
\partial_t(\vr\vu)+
\Div (\vr \vu \otimes \vu) + \Grad \pi_{\delta} (Z) = \Div \mathbb{S}(\Grad \vu) -\eta \Grad\vr\cdot\Grad\vu   +{ \eta\Div\Big(|\Grad(\vu-\vu_\infty)|^2\Grad (\vu-\vu_\infty)\Big)}
\end{equation}
\begin{equation}
\label{A4}
\vu(0,x)=\vu_0(x),\; 
\vu|_{I\times\partial \Omega} = \vu_B,
\end{equation}
\begin{equation}
\label{A5}
\partial_tZ-\eta \Del Z  + \Div (Z\vu) = 0,
\end{equation}
\begin{equation}
\label{A6}
Z(0,x)=Z_0(x),\;\left( -\eta \ \Grad Z + Z \vu \right) \cdot \vc{n}|_{I\times\partial \Omega} = \left\{
\begin{array}{l} Z_B \vu_B \cdot \vc{n} \ \mbox{if}\ [\vu_B \cdot \vc{n}](x) \leq 0,\,x\in \partial\Omega,\\
Z \vu_B \cdot \vc{n}\;\mbox{if}\ [\vu_B \cdot \vc{n}](x) > 0,\,x\in \partial\Omega,\end{array} \right.
\end{equation}
with positive parameters $\ep> 0$, $\delta > 0$, $\eta>0$.
The solution to this system is obtained by means of Galerkin approximation of the momentum equation and Banach Fixed point theorem for existence of unique local in time solutions.
We are not going to repeat all these details here, the details of this procedure can be found in \cite{CHJN} for the case with one continuity equation and it can be combined with the ideas and techniques from \cite{M3NPZ}, where two continuity equations were considered, exactly in the same setting as here: one quantity is included into the pressure, the other into the momentum. Let us only notice that when $\vu$ is replaced by it's Galerkin approximation $\vu^n$, then it is still possible to prove the comparison principle between $\vr$ and $Z$, similarly to the above mentioned paper.
Indeed, taking $c_\star,c^\star$ as in (\ref{fs-}) we may write
$$
\partial_t(Z-c_\star \vr ) - \eta \Delta_x (Z-c_\star \vr ) + \Div{\big(\vu^n(Z-c_\star \vr )\big)} = 0,
$$
and
$$
\partial_t(c^\star \vr-Z)  - \eta \Delta_x (c^\star \vr-Z) + \Div\big(\vu^n (c^\star \vr-Z)\big) = 0,
$$
with the corresponding boundary conditions. Therefore, exactly as in Lemma 4.3 from \cite{CHJN} we show that
since both equations have non-negative initial conditions, it is easy to see that also the solutions are non-negative and due to the uniqueness of solutions we deduce that
$$(Z-c_\star \vr)(t,x)\geq \inf_{x\in\overline{\Omega}}(Z_0-c_\star \vr_0)e^{-Kt},\quad (c^\star \vr-Z)(t,x)\geq \inf_{x\in\overline{\Omega}}(c^\star \vr_0-Z_0)e^{-Kt},$$
as well as $\vr(t,x) \geq \inf_{x\in \overline{ \Omega}}e^{-Kt}$, where
\bFormula{K}
K= \|{\rm div}\,\vu^n\|_{L^\infty(Q_T)}.
\eF
Using again the assumptions on the initial data 
{
\eqref{data_Z}}, we therefore obtain
$$0 < c_\star \vr \leq Z \leq c^\star \vr \mbox{ a.e. in } Q_T.$$
With these inequalities in place we can let $n\to\infty$ and $\eta\to 0$ in order to obtain the weak solutions to system \eqref{main_transformed} with $\pi_\ep$ replaced by $\pi_\delta$. Note, however, that in the limit process we do not control the $L^\infty$ norm of ${\rm div}\, \vu$ which finally leads to 
$$0 \leq c_\star \vr \leq Z \leq c^\star \vr \mbox{ a.e. in } Q_T.$$ All other steps are more or less standard except for the fact that the solution is renormalized in the sense of Definition \ref{def_renor}. This fact, however, follows directly from \cite[Lemma 3.1]{CHJN} and will return to the procedure when letting $\ep\to 0^+$.

The existence result is summarised in the following theorem.

\begin{Theorem} \label{TM1-}
Let $\Omega \subset R^d$, $d = 2,3$, be a bounded domain of class $C^{2}$ such that $\Gamma_{\rm in}$ is an open $C^2$ $d-1$ dimensional manifold. Let $\delta >0$, $\ep>0$ and $T>0$.
Under the assumptions 
{
\eqref{As1}--\eqref{Ass3}} the 
problem \eqref{ia1}--\eqref{ia4} with the pressure $\pi_\ep$ replaced by $\pi_\delta$ admits at least one renormalized bounded energy weak solution $(\vrd, \vud, \Zd)$, i.e.
\begin{description}
\item{1.} 
The triple $(\vrd,\vud,\Zd)$ belongs { to the following functional space:
\bFormula{fs}
\begin{split}
&\vrd, \Zd \in L^\infty (0,T; L^\gamma(\Omega)), \quad 0 \leq c_\star \vrd \leq \Zd \leq c^\star \vrd \ \text{ a.e. in } (0,T) \times \Omega, \\
&\vud \in L^2(0,T; W^{1,2}(\Omega; R^d)), \quad 
\vud|_{I\times\partial \Omega} = \vu_B.
\end{split}
\eF}
\item{2.}
{ The function $\vrd\in C_{\rm weak}([0,T], L^\gamma(\Omega))$  satisfies the integral identity
\bFormula{ce}
\begin{split}
&\intO{\vrd(\tau,\cdot)\varphi(\tau,\cdot)} - \intO{\vr_0(\cdot)\varphi(0,\cdot)} \\
&=
\int_0^\tau\intO{\Big(\vrd\partial_t\varphi + \vrd \vud \cdot \Grad \varphi\Big) }{\rm d}t - \int_0^\tau\int_{\Gamma_{\rm in}} \vr_B \vu_B \cdot \vc{n} \varphi 
\ {\rm d}S_x{\rm d} t 
\end{split}
\eF}
for any $\tau\in [0,T]$ and $\varphi \in C_c^1([0,T]\times({\Omega}\cup\Gamma_{\rm in}))$. 
In particular,
\bFormula{mi-}
\intO{\vrd(\tau,\cdot)}\le\intO{\vr_0}-\int_0^\tau\int_{\Gamma_{\rm in}}\vr_B\vu_B\cdot\vc n \ {\rm d}S_x{\rm d}t.
\eF
\item{\it 3.} 
{ The renormalized continuity equation 
\bFormula{rce}
\begin{split}
&\intO{b(\vrd)(\tau,\cdot)\varphi(\tau,\cdot)} - \intO{b(\vr_0)(\cdot)\varphi(0,\cdot)} \\
&=
\int_0^\tau\intO{\Big(b(\vrd)\partial_t\varphi + b(\vrd )\vud \cdot \Grad \varphi+(b(\vrd)-b'(\vrd)\vrd)\Div\vud\Big) }{\rm d}t \\
&\quad - \int_0^\tau\int_{\Gamma_{\rm in}} b(\vr_B) \vu_B \cdot \vc{n} \varphi \ {\rm d}S_x{\rm d} t \\
\end{split}
\eF
holds for any $b\in C[0,\infty)$ with $b'\in C_c[0,\infty)$, $\tau\in [0,T]$, and $\varphi \in C_c^1([0,T]\times({\Omega}\cup\Gamma_{\rm in}))$.}
\item{4.} 
{ The function $\vrd\vud\in  C_{\rm weak}([0,T], L^{\frac{2\gamma}{\gamma+1}}(\Omega;R^d))$ satisfies the integral identity
\bFormula{me}
\begin{split}
&\intO{\vrd\vud(\tau,\cdot)\cdot\bfphi(\tau,\cdot)} - \intO{\vr_0\vu_0(\cdot)\cdot\bfphi(0,\cdot)} \\
&=
\int_0^\tau\intO{  \Big(\vrd\vud\cdot\partial_t\bfphi+ \vrd \vud \otimes \vud  : \Grad \bfphi + \pi_\delta(\Zd) \Div \bfphi  
-  \mathbb{S}(\Grad \vud) : \Grad \bfphi + \vrd(\vc{w}-\vud)\cdot \bfphi\Big) }{\rm d}t 
\end{split}
\eF
for any $\tau\in [0,T]$ and $\bfphi \in C^1_c ([0,T]\times\Omega; R^d)$.}
\item{5.}
 The function $\Zd\in C_{\rm weak}([0,T], L^\gamma(\Omega))$  satisfies the integral equalities \eqref{ce}--\eqref{rce} with $\vrd$ replaced by $\Zd$, $\vr_B$ replaced by $Z_B$, and $\vr_0$ replaced by $Z_0$.
\item{6.} 
The energy inequality
\bFormula{ei}
\begin{split}
&\intO{\Big(\frac 12\vrd|\vud-\vu_\infty|^2 + { H_{\delta}}(\Zd)\Big)(\tau)} 
+ \int_0^\tau\intO{\tn S(\Grad(\vud-\vu_\infty)):\Grad(\vud-\vu_\infty)}{\rm d}t \\
&\le 
{ \intO{\Big(\frac 12\vr_0|\vu_0-\vu_\infty|^2 + { H_{\delta}}(Z_0)\Big)}}
- \int_0^\tau\intO{\pi_\delta(\Zd)\Div\vu_\infty}{\rm d}t \\
&\quad -\int_0^\tau\intO{\vrd\vud\cdot\Grad\vu_\infty\cdot(\vud-\vu_\infty)}{\rm d}t - \int_0^\tau\intO{\tn S(\Grad\vu_\infty):\Grad(\vud-\vu_\infty)}{\rm d}t \\
&\quad - \int_0^\tau\int_{\Gamma_{\rm in}} { H_{\delta}}(Z_B)\vu_B\cdot\vc n\ {\rm d}S_x{\rm d}t + \int_0^\tau\intO{\vrd(\vc{w}-\vud)\cdot(\vud-\vu_\infty)}{\rm d}t
\end{split}
\eF
holds for a.a. $\tau\in (0,T)$ and the extension $\vu_\infty\in W^{1,\infty}(\Omega;R^d)$ of $\vu_B$ discussed above.
In 
\eqref{ei}, { the function $H_{\delta}(z)$ is defined by} 
\bFormula{Hep}
H_{\delta}(z)= z\int_0^z \frac{ \pi_\delta(s)}{s^2}\ {\rm d}s.
\eF
\end{description}
\end{Theorem}

\section{Uniform estimates with respect to $\delta$ and limit $\delta\to 0$}\label{4}

In this section 
we prove our first main result, Theorem \ref{TM1!}. For this reason we will deduce some  a-priori estimates that are uniform with respect to $\delta$ with fixed positive $\ep$. Then we will improve these estimates to show that the approximation of the pressure is in fact uniformly integrable. This, together with Lions compactness argument for the density sequence will allow us to identify the limit of the pressure term, and hence the whole system.

\subsection{Uniform estimates}
\label{4.1}

We start by deriving
uniform estimates for the triple $(\vrd,\vud,\Zd)$ constructed in Theorem \ref{TM1-}.
Note that due to Lemma
\ref{Lextrem}, we have $\intO{\pi_\delta(\Zd(t,\cdot))\Div\vu_\infty}\ge 0$ at all time levels, and so, following  \cite[Section 4.3.3]{CHJN} we can show that{ the} energy inequality \eqref{ei} in combination with the conservation of mass (\ref{mi-})  yields 
\begin{equation} \label{Es3}
\begin{aligned}
\|H_{\delta}(\Zd)\|_{L^\infty(0,T;L^1(\Omega))} &\le c({\rm data}),  \\
\|\vrd|\vud|^2\|_{L^\infty(0,T; L^1(\Omega))} &\le c({\rm data}),  \\
\|\vud\|_{L^2(0,T;W^{1,2}(\Omega))} &\le c({\rm data}). 
\end{aligned}
\end{equation}
From these estimates, using \eqref{Hep}, it follows that 
\begin{equation}
\label{Es1}
\|\Zd\|_{L^\infty(0,T;L^\gamma(\Omega))} \le c({\rm data}),
\end{equation}
and hence, due to \eqref{fs}
\begin{equation}
\label{Es1r}
\|\vrd\|_{L^\infty(0,T;L^\gamma(\Omega))} \le c({\rm data}).
\end{equation}
Moreover
\begin{equation}
\label{Es1-}
{\rm ess\ sup}_{t\in(0,T)}\intO{ \tilde{H}_{\delta}(\Zd(t,x))}\le c ({\rm data}),
\end{equation}
where
\bFormula{hep}
\tilde{H}_{\delta}(z)=
{ 
\begin{cases}
\ep(1-z)^{-(\beta-1)} &\text{ if } z\in [0,1-\delta],\\
\ep\delta^{-(\beta-1)}+\ep\delta^{-\beta}(z-1+\delta) &\text{ if } z\in (1-\delta,\infty).
\end{cases}
}
\eF
This fact follows from the form of the pressure $\pi_\delta$ and the energy $H_{\delta}$, where $\tilde{H}_{\delta}(\Zd)$ contains the most singular terms in $\delta$ for $\delta \to 0+$. 
By virtue of \eqref{Es3} and \eqref{Es1r} together with \eqref{fs} 
\begin{equation} \label{Es5}
\begin{aligned}
\|\vrd\vud\|_{L^\infty(0,T; L^{\frac {2\gamma}{\gamma+1}}(\Omega))}+
\|\vrd\vud\|_{L^2(0,T; L^{\frac{6\gamma}{\gamma+6}}(\Omega))}&\le c({\rm data}), \\
\|\Zd\vud\|_{L^\infty(0,T; L^{\frac {2\gamma}{\gamma+1}}(\Omega))}+
\|\Zd\vud\|_{L^2(0,T; L^{\frac{6\gamma}{\gamma+6}}(\Omega))}&\le c({\rm data}).
\end{aligned}
\end{equation}

\subsection{Limit in the continuity equation and boundedness of density}
\label{adce}

{ We deduce from the estimates \eqref{Es3}--\eqref{Es1r} that 
\bFormula{uweak}
\begin{split}
\vu_\delta\rightharpoonup\vu \quad &\text{ in } L^2(0,T;W^{1,2}(\Omega)), \\
Z_\delta\rightharpoonup^*Z \quad &\text{ in } L^\infty(0,T; L^\gamma(\Omega)),\\
\vr_\delta\rightharpoonup^*\vr \quad &\text{ in } L^\infty(0,T; L^\gamma(\Omega)),
\end{split}
\eF
at least for a subsequence.}
We also deduce from the continuity equation \eqref{ce} and its version for $Z_\delta$,  thanks to  \eqref{Es5} and \eqref{fs},
that the sequences of functions $t \mapsto \intO{\vrd\phi}$, and $t \mapsto \intO{\Zd\phi}$, 
$\phi\in C^1_c(\Omega)$, are equi-continuous.
{ Therefore, by { the} Arzel\`a--Ascoli theorem and separability of $L^{{\gamma'}}(\Omega)$, we get 
\bFormula{wreweak}
\begin{split}
\vrd\to\vr \quad \text{ in } C_{\rm weak}(0,T;L^\gamma(\Omega)),\\
\Zd\to Z \quad \text{ in } C_{\rm weak}(0,T;L^\gamma(\Omega)).
\end{split}
\eF
Both of the sequences converge strongly in $L^2(0,T;W^{-1,2}(\Omega))$ due to compact 
embedding $L^\gamma(\Omega)\hookrightarrow\hookrightarrow 
W^{-1,2}(\Omega)$.} In particular, this implies that
\[
\vrd\vud\rightharpoonup\vr\vu \quad \text{in } L^2(0,T; L^{\frac{6\gamma}{\gamma+6}}(\Omega)),
\]
\[
\Zd\vud\rightharpoonup Z\vu \quad \text{in } L^2(0,T; L^{\frac{6\gamma}{\gamma+6}}(\Omega)).
\]

This enables us to pass to the limit in the weak formulation \eqref{ce} so that we get the identity 
\bFormula{continuity}
\begin{split}
&\intO{\vr(\tau,\cdot)\varphi(\tau,\cdot)} - \intO{\vr_0(\cdot)\varphi(0,\cdot)} \\
&=
\int_0^\tau\intO{\Big(\vr\partial_t + \vr \vu \cdot \Grad \varphi\Big) }{\rm d}t - \int_0^\tau\int_{\Gamma_{\rm in}} \vr_B \vu_B \cdot \vc{n} \varphi 
\ {\rm d}S_x{\rm d} t 
\end{split}
\eF
and its analogue for $Z$. Both hold for any $\tau\in [0,T]$ and $\varphi \in C_c^1([0,T]\times({\Omega}\cup\Gamma_{\rm in}))$.

To conclude this subsection, we deduce from \eqref{Es1-} that 
\bFormula{Es7}
0 \le Z(t,x) < 1 \quad \text{a.e. in }\Omega.
\eF
Indeed, for any fixed sufficiently small $\delta^*>0$, we have
\bFormula{nn}
\intO{\tilde H_{\delta^*}(Z(t))} \le \liminf_{\delta\to 0} \intO{\tilde H_{\delta^*}(\Zd(t))}\le \liminf_{\delta\to 0} \intO{\tilde H_{\delta}(\Zd(t))} \le  c({\rm data})
\eF
 for  almost all $t\in (0,T)$, where the first inequality  is a consequence of convexity of function $H_{\delta^*}(\cdot)$  on $[0,1-\delta)$ as well its linearity in $Z_\delta$ in the remaining part, second inequality follows from monotonicity of the map $\delta\mapsto \tilde H_\delta(Z)$  in a small right neighbourhood of $0$, and the third inequality follows from  (\ref{Es1-}). Next, as $\tilde H_\delta(\cdot)$ is globally Lipschitz, using the continuity equation we deduce that { $Z\in C([0,T);L^1(\Omega))$ and then $\tilde H_{\delta^*}(Z)\in C([0,T]; L^1(\Omega))$.} Therefore formula
(\ref{nn}) implies
\bFormula{nn+}
\intO{\tilde H_{\delta^*}(Z(t))} \le c({\rm data}) 
{ \quad \text{ for all } t\in [0,T]},
\end{equation}
and uniformly  in $\delta^*$. This implies that $Z\leq 1$. Finally letting $\delta^*\to 0$ in \eqref{nn+}, recalling \eqref{hep}, we obtain
\[
\intO{(1-Z)^{-(\beta-1)}}\le c({\rm data})
{ \quad \text{ for all } t\in [0,T],}
\]
which yields \eqref{Es7}.

\subsection{Uniform integrability of pressure}\label{4.2}

In order to  pass to the limit in the weak formulation of the momentum equation \eqref{me}, we have to improve estimates for pressure.
So far, we do not even know whether the pressure is  uniformly integrable in $\delta$. 
In this section we are going to prove it. 

A general tool to obtain these estimates is the following Bogovskii lemma (see, e.g., \cite{GALN} or \cite[Theorem 10.11]{FEINOV}).

{ 
\bLemma{Bog}
Let $\Omega \subset R^d$, $d\geq 2$, be a bounded Lipschitz domain. Then there exists a linear operator
\[
{\cal B}:\Big\{f\in C^\infty_c(\Omega;R^d)\,\colon\,\intO{f}=0\Big\} \mapsto C^\infty_c(\Omega;R^d) 
\]
satisfying the following three properties.
\begin{enumerate}
\item
For all $f\in C^\infty_c(\Omega;R^d)$ satisfying $\intO{f}=0$
\[
{\rm div}{\cal B}[f]=f.
\]
\item
Let $\overline L^p(\Omega):=\{f\in L^p(\Omega)\,|\,\intO{f}=0\}$.
Then the operator ${\cal B}$ extends to a bounded linear operator from $\overline{L}^p(\Omega)$ to $W^{1,p}(\Omega)$ for any $1<p<\infty$. 
In other words, for each $1<p<\infty$ there is $c(p)>0$ such that for all $f\in \overline L^p(\Omega)$
\[
\|{\cal B}[f]\|_{W^{1,p}(\Omega;R^3)}\le c(p)\|f\|_{L^p(\Omega)}.
\]
\item
If $f={\rm div}\,\vc g$ for some $\vc g\in L^q(\Omega)$, $1<q<\infty$ with $\vc g\cdot\vc n|_{\partial\Omega}=0$ in the sense of normal traces, then
there is $c(q)>0$ such that 
\[
\|{\cal B}[f]\|_{L^q(\Omega;R^3)}\le c(q)\|\vc g\|_{L^q(\Omega,R^3)}
\] 
for all $\vc g$ with the above properties.
\end{enumerate}
\eL 
}

We employ this lemma to construct suitable test functions for the momentum equation. Note that by standard density argument we can extend the class of test functions in \eqref{me} to certain $W^{1,q}$-functions with zero trace in $\Gamma_{\rm out}$. Our test function will be a suitable test function due to estimates performed below.

We use in \eqref{me} the following test function
\begin{equation}\label{B1}
\vcg{\varphi} = \eta(t) {\cal B}\Big(\psi Z_\delta-\frac{\psi}{\intO{\psi}} \intO{\psi Z_\delta}\Big),    
\end{equation}
where $\eta \in C^1_c(0,T)$ and $\psi \in C^\infty_c(\Omega)$, $0\leq \eta,\psi\leq 1$. Then we have
\begin{equation} \label{B2}
\int_0^T \intO{\eta \pi_\delta(Z_\delta) \Big(\psi Z_\delta - \frac{\psi}{\intO{\psi}}\intO{\psi Z_\delta}\Big)}{\rm d}t = \sum_{i=1}^6 I_i,
\end{equation}
where
\begin{align*}
I_1&=-\int_0^T\partial_t\eta\intO{\vrd\vud\cdot{\cal B}\Big(\psi\Zd - \frac{\psi}{\intO{\psi}}\intO{\psi Z_\delta}\Big)}{\rm d}t, \\
I_2&=\int_0^T\eta\intO{\vrd\vud\cdot{\cal B}({\rm div}\,(\Zd\vud\psi))}{\rm d}t, \\
I_3&=-\int_0^T\eta\intO{\vrd\vud\cdot{\cal B}\Big(\Zd\vud\cdot\Grad\psi
- \frac{\psi}{\intO{\psi}}\intO{\Zd\vud\cdot\Grad\psi}\Big)}{\rm d}t, \\
I_4&=-\int_0^T\eta\intO{\vrd(\vud\otimes\vud):\Grad{\cal B}\Big(\psi\Zd - \frac{\psi}{\intO{\psi}}\intO{\psi Z_\delta}\Big)}{\rm d}t, \\
I_5&=\int_0^T\eta\intO{\tn S(\Grad\vud):
\Grad{\cal B}\Big(\psi\Zd - \frac{\psi}{\intO{\psi}}\intO{\psi Z_\delta}\Big)}{\rm d}t, \\
I_6&= \int_0^T\eta\intO{\vrd (\vc{w}-\vud)\cdot
{\cal B}\Big(\psi\Zd - \frac{\psi}{\intO{\psi}}\intO{\psi Z_\delta}\Big)}{\rm d}t.
\end{align*}

Clearly, $|\sum_{i=1}^6 I_i| \leq C$ with $C$ independent of $\delta$ by estimates \eqref{Es3}--\eqref{Es1r} and we end up with
\begin{equation} \label{B3}
\int_0^T \intO{\eta \pi_\delta(Z_\delta) \Big(\psi Z_\delta - \frac{\psi}{\intO{\psi}}\intO{\psi Z_\delta}\Big)}{\rm d}t \leq C.
\end{equation}
We now choose $\psi \in C^\infty_c(\Omega)$ such that  $0\leq \psi\leq 1$ and $\psi \equiv 1$ in $K$ some compact set $K\subset\subset \Omega$.
For the fixed set $K$ and any such $\psi$ as above we further denote
$$
M_{\delta,K}=\max_{t\in [0,T]} \intO{Z_\delta \psi}.
$$
We claim that for each $K \subset\subset \Omega$ there exists $\delta_0 >0$ and $\lambda >1$ such that for any $\delta <\delta_0$ it holds
\begin{equation}\label{lambda}
\lambda M_{\delta,K} < \intO{\psi}.
\end{equation}
This follows from several facts shown before. The form of  $\tilde H_\delta(Z_\delta)$, see \eqref{hep}, yields that for any positive (sufficiently small) $\delta_\theta$ we have on the set $Z_\delta>1-\delta_\theta$ (recall that $\ep$ is fixed at this moment)
\eq{
\tilde H_\delta(Z_\delta)> \frac{C}{\delta_\theta^{\beta-1}}.
}
Since $\beta>\frac 52$ (if $d=3$; note that it is enough to have $\beta >2$ which is the case if $d=2$) and due to the $L^\infty(0,T;L^1(\Omega))$ bound of $\tilde H_\delta(Z_\delta)$ we see that for arbitrarily small $\theta>0$ there exists $\delta_\theta$ such that
$$
\sup_{\delta<\delta_\theta} \sup_{t\in[0,T]}\left|x\in \Omega: Z_\delta(t,x) > 1-\delta_\theta\right|<\theta.
$$
Furthermore, for $\theta>0$, sufficiently small, the $L^\infty(0,T;L^1(\Omega))$ bound of $\tilde H_\delta(Z_\delta)$ implies that we may take $\delta_\theta = \theta^{\frac 23}$. 
Therefore we have (we assume $\delta \leq \delta_\theta$):
$$
\begin{aligned}
\max_{t\in[0,T]}\intO{Z_\delta \psi}&=\sup_{t\in[0,T]}\int_{\{Z_\delta(t,x) > 1-\delta_\theta\}}Z_\delta \psi \dx+\sup_{t\in[0,T]}\int_{\{Z_\delta(t,x) \leq 1-\delta_\theta\}}Z_\delta \psi \dx\\
&\leq \theta^{\frac{1}{\gamma'}}\|Z_\delta\|_{L^\infty(0,T;L^\gamma(\Omega))}+(1-\delta_\theta)\intO{\psi}\\
&\leq C\theta^{\frac{1}{\gamma'}}+(1-\delta_\theta)\intO{\psi}.
\end{aligned}
$$
Thus, taking $\theta$ possibly even smaller, we may achieve by taking $\gamma$ sufficiently large that $C\theta^{\frac{1}{\gamma'}} <\frac {\delta_\theta}{2} \intO{\psi}$  which leads to
\begin{equation} \label{est_a}
\max_{t\in[0,T]}\intO{Z_\delta \psi}\leq\lr{1-\frac{\delta_\theta}{2}}\intO{\psi}.
\end{equation}
Thus, for $\lambda:= \frac{1}{1-\frac{\delta_\theta}{2}}$ we showed \eqref{lambda}.


We return to inequality \eqref{B3}. Clearly, for the set $O_1:= \Big\{(t,x) \in (0,T)\times \Omega: Z_\delta(t,x) < \frac{\lambda M_{\delta,K}}{\intO{\psi}}<1\Big\}$ we have
$$
\begin{aligned}
\Big|\int\int_{O_1} \eta \pi_\delta(Z_\delta) &\Big(\psi Z_\delta-\frac{\psi}{\intO{\psi}} \intO{Z_\delta\psi}\Big)\ {\rm d}x{\rm d}t\Big| \\
&\leq \Big|\int\int_{O_1} \eta \pi_\delta\Big(\frac{\lambda M_{\delta,K}}{\intO{\psi}}\Big) \Big(\psi Z_\delta-\frac{\psi}{\intO{\psi}} \intO{Z_\delta\psi}\Big)\ {\rm d}x{\rm d}t\Big| \leq CT.
\end{aligned}
$$
Next
$$
\begin{aligned}
&\int\int_{((0,T)\times \Omega)\setminus O_1} \eta\pi_\delta(Z_\delta) \psi Z_\delta \ {\rm d}x{\rm d}t \leq C +  \int\int_{((0,T)\times \Omega)\setminus O_1}\eta \pi_\delta(Z_\delta) \Big(\frac{\psi}{\intO{\psi}} \intO{\psi Z_\delta}\Big) \ {\rm d}x{\rm d}t \\ 
&\leq C + \int\int_{((0,T)\times \Omega)\setminus O_1} \eta\psi\pi_\delta(Z_\delta)\frac{M_\delta}{\intO{\psi}} \ {\rm d}x{\rm d}t \leq C + \frac{1}{\lambda} \int\int_{((0,T)\times \Omega)\setminus O_1} \eta\psi \pi_\delta(Z_\delta) Z_\delta \ {\rm d}x{\rm d}t.
\end{aligned}
$$
The computations above imply that for any $K \subset\subset \Omega$ and corresponding $\psi\in C^\infty_c(\Omega)$ there exists $C=C(K)$ such that
\begin{equation}\label{B4}
\int_0^T \eta \intO{\psi \pi_\delta(Z_\delta) Z_\delta}{\rm d}t \leq C(K).
\end{equation}
Whence we also have
\begin{equation}\label{B5}
\int \int_{\{Z_\delta \leq 1-\delta \}} \eta \psi \ep (1-Z_\delta)^{-\beta} \ {\rm d}x {\rm d}t \leq C(K),
\end{equation}
where $\eta \in C^\infty_c(0,T)$.
\begin{Remark}\label{RB1}
Note that this estimate was obtained without assuming short time interval for zero velocity flux as it was the case in \cite{CHNY2}. Therefore, from this point of view our paper even improves the result in the above cited paper in the sense that global in time solution exists provided 
$$
\int_{\partial \Omega} \vu_B \cdot \vc{n}\ {\rm d}S\geq 0.
$$
The case of negative flux, however, remains an interesting open problem.
\end{Remark}

\subsection{Equi-integrability of pressure}
\label{4.4}

In order to show equi-integrability of the sequence $\pi_\delta(\Zd)$, we shall use the renormalized continuity equation \eqref{rce} with $\vrd$ replaced by $\Zd$.
We fix the same cut-off functions $\eta$ in the time variable and $\psi$ in the spatial variables as in the previous section and
$0\le \psi\in C^1_c(\Omega)$ and consider the following test function
\[
\vcg{\varphi}=\eta(t){\cal B}(\psi b(\Zd)-\alpha_\delta) \quad \text{ where } \alpha_\delta = \frac 1{|\Omega|}\intO{\psi b(\Zd)}
\]
with 
\begin{align*}
b(Z) 
&= \begin{cases}
-\ln(1/2) & \text{ if } Z\in [0,1/2], \\
-\ln(1-Z) & \text{ if } Z \in (1/2,1-\delta), \\
-\ln\delta & \text{ if }Z \in [1-\delta,\infty).
\end{cases} 
\end{align*}
We note that 
\[
b'(Z) = \frac 1{1-Z}1_{(1/2,1-\delta)}(Z),
\]
where, as above, $1_E(Z)$ denotes the characteristic function of a set $E$.
In view of  \eqref{Es1}, \eqref{Es1-}, and \eqref{B5}, we notice also that for any $1\le p<\infty$,  and any compact $K\subset\Omega$,
\bFormula{pr3}
\begin{split}
\|b(\Zd)\|_{L^ \infty(0,T; L^p(K))} &\le c({\rm data}, K, p), \\
\|\eta \Zd b'(\Zd)-b(\Zd)\|_{L^\beta((0,T)\times K)} &\le c({\rm data},K,\eta),\\
\|\Zd b'(\Zd)-b(\Zd)\|_{L^\infty(0,T; L^{\beta-1}(K))} &\le c({\rm data},K),
\end{split}
\eF
where $\eta$ is as above.
We test  the momentum equation \eqref{me}  by $\vcg{\varphi}$ to obtain the following identity
\[
\int_0^T \eta \intO{\psi \ep\pi_\delta(\Zd) b(\Zd)}{\rm d}t= \sum_{i=1}^{8}I_i,
\]
where
\begin{align*}
&I_1=\frac 1{|\Omega|}\int_0^T\eta(t)\int_{\Omega}\psi b(\Zd){\rm d}x\intO{\ep\pi_\delta(\Zd)}{\rm d}t, \\
&I_2=-\int_0^T\partial_t\eta\intO{\vre\vue\cdot{\cal B}(\psi b(\Zd)-\alpha_\delta)}{\rm d}t, \\
&I_3=\int_0^T\eta\intO{\vrd\vud\cdot{\cal B}\Big(\Div(\psi b(\Zd)\vud)\Big)}{\rm d}t, \\
&I_4=-\int_0^T\eta\intO{\vrd\vud\cdot{\cal B}\Big( b(\Zd)\vud\cdot\Grad\psi- \frac 1{|\Omega|}\int_{\Omega} b(\Zd)\vud\cdot\Grad\psi \ {\rm d}x\Big)    }{\rm d}t, \\
&I_5=\int_0^T\eta\intO{\vrd\vud\cdot{\cal B}\Big[\psi\Big(\Zd b'(\Zd)-b(\Zd)\Big)\Div\vud
- \frac 1{|\Omega|}\intO{\psi\Big(\Zd b'(\Zd)-b(\Zd)\Big)\Div\vud}\Big]}{\rm d}t, \\
&I_6=\int_0^T\eta\intO{\tn S(\Grad\vud):\Grad{\cal B}(\psi b(\Zd)-\alpha_\delta)}{\rm d}t, \\
&I_7=- \int_0^T\eta\intO{\vrd\vud\otimes\vud:\Grad{\cal B}(\psi b(\Zd)-\alpha_\delta)}{\rm d}t, \\
&I_8 = -\int_0^T\eta \intO{\vrd (\vc{w}-\vud)}\cdot {\cal B}(\psi b(\Zd)-\alpha_\delta){\rm d}t.
\end{align*}

The above calculation involves integration by parts and the renormalized equation \eqref{rce} for unknown $Z_\delta$.
The function $b$ is clearly admissible in the renormalized continuity equation.
We verify, using the approach of \cite{FEMAL}, estimates  \eqref{Es3}--\eqref{Es1-}, and \eqref{pr3}, that {for any $\beta>5/2$ there is $\gamma>3/2$ (sufficiently large - $\gamma\to\infty$ as $\beta\to \frac 52+$)} such that absolute values of $I_1,\dots,I_8$ are bounded
above by some positive constants.
The most severe constraints on the values of $\beta$
and $\gamma$ within these calculations are imposed  in estimating the term $|I_5|$.  Note that in case $d=2$ the same approach as in \cite{FEMAL} can be used to see that the most singular term can be estimated for $\beta >2$. 
Effectuating this process, we obtain that for any compact set $K\subset \Omega$,
\bFormula{pvre}
\|\eta\pi_\delta(\Zd)b(\Zd)\|_{L^{1}((0,T)\times K)}\le c({\rm data},K, \eta).
\eF
Consequently,  for fixed $\ep>0$, the sequence $\pi_\delta(\Zd)$ is equi-integrable in $L^{1}(Q)$ for any $Q\subset\subset (0,T)\times\Omega$ and
\eq{\label{pvrecon}
\pi_\delta(\Zd)\rightharpoonup\overline{\pi(Z)} \quad \text{ in } L^{1}(J\times K) 
}
for any compact set $J\times K\subset (0,T)\times \Omega$ at least for a chosen subsequence (not relabeled).

\subsection{Momentum equation}
\label{4.5}

With the help of \eqref{Es3}, \eqref{Es1r}, and \eqref{pvre} employed in the momentum equation \eqref{me}, we verify equicontinuity of the sequence $t\mapsto \intO{\vrd\vud(t,\cdot)\varphi(\cdot)}$ in $C[0,T]$ for any $\varphi\in C^1_c(\Omega)$ .
Therefore, we may use the Arzel\`a--Ascoli theorem in combination with \eqref{Es5} the separability of $L^{\lr{\frac{2\gamma}{\gamma+1}}'}(\Omega)$ to show that 
\bFormula{bweak}
\vrd\vud\to\vr\vu \quad \text{ in } C_{\rm weak}([0,T];L^{\frac{2\gamma}{\gamma+1}}(\Omega)).
\eF
Consequently, the imbedding $L^{\frac{2\gamma}{\gamma+1}}(\Omega)\hookrightarrow\hookrightarrow W^{-1,2}(\Omega)$ (for $\gamma>\frac32$) in combination with the weak convergence
of $\vu_n$ in $L^2(0,T; W^{1,2}(\Omega))$ { implies}
\bFormula{cweak}
\vrd\vud\otimes\vud\rightharpoonup\vr\vu\otimes\vu \quad \text{ in } L^2(0,T; L^{{6\gamma}/{4\gamma+3}}(\Omega)).
\eF

Thus, letting $\delta\to 0$ in weak formulation of \eqref{me} while using 
\eqref{uweak}, \eqref{wreweak}, \eqref{pvrecon}, and \eqref{cweak} we obtain that for any $\tau\in [0,T]$ and $\bfphi \in C^1_c ([0,T]\times\Omega; R^d)$,
\bFormula{me+}
\begin{split}
&\intO{\vr\vu(\tau,\cdot)\cdot\bfphi(\tau,\cdot)} - \intO{(\vr_0\vu_0)(\cdot)\cdot \bfphi(0,\cdot)} 
= \intO{  \left( \vr \vu \otimes \vu  : \Grad \bfphi + \ep\overline{\pi(Z)} \Div \bfphi\right)} \\ - &\int_0^\tau\intO{ \mathbb{S}(\Grad \vu) : \Grad \bfphi }{\rm d}t + \int_0^\tau \intO{\vr (\vc{w}-\vu)\cdot \bfphi}{\rm d}t.
\end{split}
\eF

The proof of Theorem \ref{TM1!} is therefore complete if we show that 
\bFormula{pp}
\overline{\pi(Z)}=\pi_\ep(Z),
\eF
which amounts, in fact, to show that the  sequence $\Zd$ converges almost everywhere in $Q_T$.

\subsection{Strong convergence of $\Zd$}\label{4.6}

We denote by $\Grad\Delta^{-1}$ the pseudodifferential operator of the 
Fourier symbol $\frac {{\rm i}\xi}{|\xi|^2}$ and by ${\cal R}$ the Riesz transform of the 
Fourier symbol $\frac {\xi\otimes\xi}{|\xi|^2}$.
Following  Lions \cite{LI4} with modified  in  \cite{FNP}, we shall use the test function
\[
\varphi(t,x)=\eta(t)\psi(x)\Grad\Delta^{-1}(\Zd\psi),\;\;\eta\in C^1_c(0,T),\;\psi\in C^1_c(\Omega)
\]
in the approximate momentum equation \eqref{me} and the test function 
\[
\varphi(t,x)=\eta(t)\psi(x)\Grad\Delta^{-1}(Z\psi),\;\;\eta\in C^1_c(0,T),\;\psi\in C^1_c(\Omega)
\]
in the limit momentum equation \eqref{me+}, subtract the resulting identities, and then perform the limit $\delta\to 0$. 
{ These calculations are laborious but nowadays standard. One can find details e.g. in  \cite[Lemma 3.2]{FNP}, \cite{NOST}, \cite{EF70} or \cite[Chapter 3]{FEINOV}) to obtain the following identity }
\bFormula{ddd!}
\begin{split}
&\int_0^T\intO{\eta\psi^2\Big(\overline{\pi(Z)} \;Z-(2\mu +\lambda) \Div\vu\, Z\Big)}\,{\rm d}t \\
&\quad -\int_0^T\intO{\eta\psi^2\Big(\overline{\pi(Z)\,Z}
-(2\mu +\lambda)\overline{Z\,\Div\vu}\Big)}{\rm d}t \\
&= \int_0^T\eta\intO{\psi^2\vu\cdot\Big(Z {\cal R}\cdot(\vr\vu)-\vr\vu\cdot{\cal R}(Z)\Big)}{\rm d}t \\
&\quad - \lim_{\delta\to 0}\int_0^T\eta\intO{\psi^2\vud\cdot\Big( \Zd {\cal R}\cdot(\vrd\vud)-\vrd\vud\cdot{\cal R}(\Zd)\Big)
}{\rm d}t.
\end{split}
\eF
In (\ref{ddd!}) and in the sequel the overlined quantities $\overline{b(Z,\vu)}$, resp. $\overline{b(Z)}$ denote $L^1(Q_T)$-weak limits of sequences
$b(\Zd,\vud)$ resp. $b(\Zd)$ (or $b_\delta(\Zd)$ if this is the case).

The most non-trivial 
moment in this process is to show that the right-hand side of this identity vanishes.
The details of this calculation and reasoning can be found in 
\cite[Lemma 3.2]{FNP}, \cite{EF70}, \cite{NOST}, or \cite[Chapter 3]{FEINOV}. Consequently,
\bFormula{ad1}
\lr{\overline{\pi(Z)\,Z} -\overline{\pi(Z)}\;Z}=(2\mu +\lambda)\Big(\overline{Z\,\Div\vu}-Z\Div\vu\Big),
\eF
and so
\bFormula{evf+}
(2\mu +\lambda) \int_0^\tau\intO{\Big(Z{\rm div}\vu-\overline{Z {\rm div}\vu}\Big)}{\rm d}t
= \int_0^\tau\intO{\Big(\overline{ \pi(Z)Z}-\overline{\pi(Z)}Z\Big)}{\rm d}t\leq 0,
\eF
due to the fact that $\pi(Z)$ is increasing.

The next (and the last) step in the proof follows closely Section 4.5 in \cite{CHJN}. Note, that this procedure does not depend on the momentum equation anymore, therefore  it will hold also for the limit passage $\ep\to0$.

Since both $(\Zd,\vud)$ and $(Z,\vu)$ satisfy the renormalized continuity equation \eqref{rce}, we obtain, in particular, that 
\begin{align*}
&\int_{\Omega}\overline{L(Z(\tau))}\varphi \ {\rm d}x-\int_{\Omega}L(Z_0)\varphi(0,x) \ {\rm d}x \\
&= \int_0^\tau\int_{\Omega} \Big(\overline{L(Z)}\partial_t\varphi+\overline{L(Z)} \vu \cdot \Grad \varphi-\varphi \overline{Z \Div \vu} \Big)\ \dx{\rm d}t
+\int_0^\tau\int_{\Gamma_{\rm in}}L(Z_B)\vu_B\cdot\vc n\varphi \ {\rm d}S_x{\rm d}t,
\end{align*}
and
\begin{align*}
&\int_{\Omega}L(Z(\tau,x))\varphi(\tau,x) \ {\rm d}x-\int_{\Omega}L(Z_0)\varphi(0,x)\ {\rm d}x \\
&= \int_0^\tau\int_{\Omega} \Big(L(Z)\partial_t\varphi+L(Z) \vu \cdot \Grad \varphi-\varphi Z\Div \vu \Big) \dx{\rm d}t
+\int_0^\tau\int_{\Gamma_{\rm in}}L(Z_B)\vu_B\cdot\vc n\varphi \ {\rm d}S_x{\rm d}t
\end{align*}
where $L(Z)=Z\ln Z$,
and $\varphi\in C^1_c([0,T]\times(\Omega\cup\Gamma_{\rm in}))$. 
Subtracting these inequalities, we obtain
\begin{align*}
&\int_{\Omega}\Big(\overline{L(Z)}-L(Z)\Big)(\tau)\varphi(\tau,x) \ {\rm d}x \\
&= \int_0^\tau\int_{\Omega} \Big(\overline{L(Z)} -L(Z)\Big)(\vu \cdot \Grad \varphi+\partial_t\varphi) \ {\rm d}x{\rm d}t-\int_0^T\intO{\varphi \Big(\overline{Z \Div \vu} 
-Z{\rm div}\vu\Big)}{\rm d}t.
\end{align*}
Hence, by virtue of \eqref{evf+} and using also in particular the function $\varphi$ independent of time, we get 
\bFormula{dod5}
\begin{split}
&\intO{\Big(\overline{Z\log Z}- Z\log Z
{ \Big)}
(\tau,x)\varphi(x)}
+ \int_0^\tau\intO{ \Big(Z\log Z-\overline{Z \log Z}\Big)  \vu \cdot \Grad \varphi }{\rm d}t\leq 0
\end{split}
\eF
for any $\tau\in [0,T]$ and $\varphi \in C^1_c(\Omega\cup\Gamma_{\rm in})$ with $\varphi \geq 0$. To show that the above formula for $\varphi(x)\to 1$ gives
\bFormula{dod6}
\intO{\Big(\overline{Z\log Z}- Z\log Z
{ \Big)}
(\tau,\cdot)}\le 0,
\eF
we can follow step by step the procedure described in Section 4.7 in \cite{CHNY2}.
On the other hand, we have
\[
\overline{Z\log Z}- Z\log Z\ge 0 \quad \text{ a.e. in } Q_T
\]
since $Z\log Z$ is convex.
Thus, 
(\ref{dod6}) yields
\[
\overline{Z\log Z}= Z\log Z \quad \text{ a.e. in } Q_T,
\]
and so
\bFormula{dod7}
\Zd\to Z \quad \text{ a.e. in } Q_T \text{ and in } L^p(Q_T) \text{ for } 1\le p<\infty,
\eF
cf. e.g. \cite[Theorem 10.20]{FEINOV}.
We deduce from \eqref{dod7} and \eqref{pvre} that for any compact $K\subset \Omega$,
\bFormula{cp}
\pi_\delta(\Zd)\to \pi(Z) \;\text{a.e. in $Q_T$ and in $L^1((0,T)\times K)$}.
\eF
In particular, we have $\overline{\pi(Z)}=\pi(Z)$ in equation \eqref{me+}, note, however, that this information is restricted solely to $Z$ and does not imply the strong convergence of $\vrd$.

\subsection{Energy inequality}

We first integrate \eqref{ei} over $0<\tau_1<\tau_2<T$ to obtain that
\bFormula{ei!}
\begin{split}
&\int_{\tau_1}^{\tau_2}\int_{\Omega}\Big(\frac 12\vrd|\vud-\vu_\infty|^2+H_\delta(\Zd)\Big)(\tau,\cdot) \ {\rm d} x{\rm d}\tau+
\int_{\tau_1}^{\tau_2}\int_0^\tau\int_{\Omega}\tn S(\Grad(\vud-\vu_\infty)):\Grad(\vud-\vu_\infty) \ {\rm d} x{\rm d}t{\rm d}\tau \\
&\le 
{ \int_{\tau_1}^{\tau_2}\intO{\Big(\frac 12\vr_0|\vu_0-\vu_\infty|^2+H_\delta(Z_0)\Big)}}{\rm d}\tau
-\int_{\tau_1}^{\tau_2}\int_0^\tau\int_{\Omega} \pi_\delta(Z_\delta)\Div\vu_\infty \ {\rm d} x{\rm d}t{\rm d}\tau \\
&\quad -\int_{\tau_1}^{\tau_2}\int_0^\tau\int_{\Omega}\vrd\vud\cdot\Grad\vu_\infty\cdot(\vud-\vu_\infty) \ {\rm d} x{\rm d}t{\rm d}\tau
-\int_{\tau_1}^{\tau_2}\int_0^\tau\int_{\Omega}\tn S(\Grad\vu_\infty):\Grad(\vud-\vu_\infty) \ {\rm d} x{\rm d}t{\rm d}\tau \\
&\quad - \int_{\tau_1}^{\tau_2}\int_0^\tau\int_{\Gamma_{{\rm in}}} H_\delta(Z_{B})\vu_{B}\cdot\vc n \ {\rm d}S_x{\rm d}t{\rm d}\tau + \int_{\tau_1}^{\tau_2}\int_0^\tau\int_{\Omega}\vrd (\vc{w}-\vud)\cdot (\vud-\vu_\infty) \ {\rm d} x{\rm d}t{\rm d}\tau.
\end{split}
\eF
Now, we can use the convergences established in Section \ref{adce} and in \eqref{dod7} and (\ref{cp}) at the right-hand side and the same convergences in combination with the lower weak semi-continuity of convex functionals
at the left-hand side (see e.g. \cite[Theorem 10.20]{FEINOV}). Note, however, that we must be slightly careful with the pressure term, where the convergence is only local. We therefore get
\bFormula{ei!++}
\begin{split}
&\int_{\tau_1}^{\tau_2}\int_{\Omega}\Big(\frac 12\vr|\vu-\vu_\infty|^2+H_\ep(Z)\Big)(\tau,\cdot) \ {\rm d} x{\rm d}\tau+
\int_{\tau_1}^{\tau_2}\int_0^\tau\int_{\Omega}\tn S(\Grad(\vu-\vu_\infty)):\Grad(\vu-\vu_\infty)\ {\rm d} x{\rm d}t{\rm d}\tau \\
&\quad  +\int_{\tau_1}^{\tau_2}\int_\alpha^{\tau-\alpha}\int_{K} \pi_\ep(Z)\Div\vu_\infty \ {\rm d} x{\rm d}t{\rm d}\tau \le 
{ \int_{\tau_1}^{\tau_2}\intO{\Big(\frac 12\vr_0|\vu_0-\vu_\infty|^2+H_\ep(Z_0)\Big)}}{\rm d}\tau
 \\
&\quad -\int_{\tau_1}^{\tau_2}\int_0^\tau\int_{\Omega}\vr\vu\cdot\Grad\vu_\infty\cdot(\vu-\vu_\infty) \ {\rm d} x{\rm d}t{\rm d}\tau
-\int_{\tau_1}^{\tau_2}\int_0^\tau\int_{\Omega}\tn S(\Grad\vu_\infty):\Grad(\vu-\vu_\infty) \ {\rm d} x{\rm d}t{\rm d}\tau \\
&\quad - 
\int_{\tau_1}^{\tau_2}\int_0^\tau\int_{\Gamma_{{\rm in}}} H_\ep(Z_{B})\vu_{B}\cdot\vc n \ {\rm d}S_x{\rm d}t{\rm d}\tau + \int_{\tau_1}^{\tau_2}\int_0^\tau\int_{\Omega}\vr (\vc{w}-\vu)\cdot (\vu-\vu_\infty) \ {\rm d} x{\rm d}t{\rm d}\tau.
\end{split}
\eF
Since the inequality holds for any $\alpha>0$, sufficiently small, and any $K$ compact subset of $\Omega$, we easily obtain at the end
\bFormula{ei!+}
\begin{split}
&\int_{\tau_1}^{\tau_2}\int_{\Omega}\Big(\frac 12\vr|\vu-\vu_\infty|^2+H_\ep(Z)\Big)(\tau,\cdot) \ {\rm d} x{\rm d}\tau+
\int_{\tau_1}^{\tau_2}\int_0^\tau\int_{\Omega}\tn S(\Grad(\vu-\vu_\infty)):\Grad(\vu-\vu_\infty)\ {\rm d} x{\rm d}t{\rm d}\tau \\
&\le 
{ \int_{\tau_1}^{\tau_2}\intO{\Big(\frac 12\vr_0|\vu_0-\vu_\infty|^2+H_\ep(Z_0)\Big)}}{\rm d}\tau
-\int_{\tau_1}^{\tau_2}\int_0^\tau\int_{\Omega} \pi_\ep(Z)\Div\vu_\infty \ {\rm d} x{\rm d}t{\rm d}\tau \\
&\quad -\int_{\tau_1}^{\tau_2}\int_0^\tau\int_{\Omega}\vr\vu\cdot\Grad\vu_\infty\cdot(\vu-\vu_\infty) \ {\rm d} x{\rm d}t{\rm d}\tau
-\int_{\tau_1}^{\tau_2}\int_0^\tau\int_{\Omega}\tn S(\Grad\vu_\infty):\Grad(\vu-\vu_\infty) \ {\rm d} x{\rm d}t{\rm d}\tau \\
&\quad - 
\int_{\tau_1}^{\tau_2}\int_0^\tau\int_{\Gamma_{{\rm in}}} H_\ep(Z_{B})\vu_{B}\cdot\vc n \ {\rm d}S_x{\rm d}t{\rm d}\tau + \int_{\tau_1}^{\tau_2}\int_0^\tau\int_{\Omega}\vr (\vc{w}-\vu)\cdot (\vu-\vu_\infty) \ {\rm d} x{\rm d}t{\rm d}\tau.
\end{split}
\eF

\subsection{Renormalized continuity equation}\label{ren_con_eq_delta}

In this section following \cite{CHJN} we can generalize the DiPerna-Lions theory for continuity equation with nonhomogenous boundary data. Recall that due to \eqref{uweak} $\vr \in L^2(0,T;L^\gamma(\Omega))$ and $\vu \in L^2(0,T;W^{1,2}(\Omega))$, and in particular $\gamma>2$. Due to \cite[Lemma~3.1]{CHJN} we may formulate the following result:
\begin{Lemma}\label{Renorm}
Let $\Omega \subset R^d$, $d=2,3$, be bounded domain of class $C^2$ such that $\Gamma_{\rm in}$ is an open $C^2$ $d-1$ dimensional manifold. 
Let $\vr_B$ and $\vu_B$ satisfy assumptions \eqref{Ass1}. Let $\vr \in L^2(0,T;L^\gamma(\Omega))$ with $\gamma >2$ and $\vu \in L^2(0,T;W^{1,2}(\Omega))$ 
 satisfy the continuity equation in a weak sense, as in \eqref{continuity}.

Then $(\vr, \vu)$ is also a renormalized solution of the continuity equation \eqref{continuity}, namely it verifies \eqref{P3} with $\vr$ instead of $Z$.
\end{Lemma}
For the detailed proof see \cite[Section~3.1]{CHJN}. Let us present here just a sketch of it. The main difficulty is to reconstruct the the boundary term on $\Gamma_{\rm in}$. To this end the set $\Omega\cup \Gamma_{\rm in}$ is extended by properly chosen open (nonempty) set $\tilde{U}^+_h(\Gamma_{\rm in})$ being in neighborhood of $\Gamma_{\rm in}$ (one can think about it as an thin pillow attached to $\Gamma_{\rm in}$ outside of $\Omega$). It is defined as follows:
\begin{equation}\label{setU}
\tilde{U}^+_h(\Gamma_{\rm in}) = 
\{ x \in {U}^+_h(\Gamma_{\rm in}) \,|\, x=\vc{X}(s,\vc{x}_0) \mbox{ for a certain }\vc{x}_0 \in\Gamma_{\rm in } \mbox{ and } 0<s<h \}
    \end{equation}
where 
$${U}^+_h(\Gamma_{\rm in}) := \{ \vc{x}_0 + z \vc{n}(\vc{x_0}) \, | \, 0<z<h, \vc{x}_0 \in\Gamma_in\} \cap (\R^d \setminus \Omega)$$
and 
$$\vc{X}'(s,\vc{x}_0) = - \tilde\vu_B (\vc{X}(s,\vc{x}_0)),\ \vc{X}(0) = \vc{x}_0 \in {U}^+_h(\Gamma_{\rm in})\cup\Gamma_{\rm in} 
\mbox{ for } s>0, \  \vc{X}(s,\vc{x}_0) \in {U}^+_h(\Gamma_{\rm in})
$$
with $$\tilde\vu_B(x) = \vu_B(\vc{x}_0), \ x = \vc{x}_0 + z \vc{n}(\vc{x}_0) \in {U}^+_h(\Gamma_{\rm in})$$
The set $\tilde{U}^+_h(\Gamma_{\rm in})$ is nonempty and open, for details see \cite[Section~3.1]{CHJN}.
Here the regularity of $\Gamma_{\rm in }$ is used. Moreover proper extension $\tilde\vu_B$ and $\tilde\vr_B$ of $\vu_B$ and $\vr_B$ on $\tilde{U}^+_h(\Gamma_{\rm in})$ is constructed such that $\tilde\vu_B \in C^1(\tilde{U}^+_h(\overline{\Gamma_{\rm in})})$, $\tilde\vr_B \in W^{1,\infty}(\tilde{U}^+_h({\Gamma_{\rm in}}))$, and 
    \begin{equation}\label{ggg}
    \Div (\tilde\vr_B \tilde\vu_B ) = 0 \mbox{ in } \tilde{U}^+_h(\overline{\Gamma_{\rm in})}, \quad \tilde\vr_B|_{\Gamma_{\rm in}} = \vr_B, \  \tilde\vu_B|_{\Gamma_{\rm in}} = \vu_B,
    \end{equation}
Then extension of $(\vr,\vu)$ to the set $\Omega_h:=\Omega\cup\Gamma_{\rm in} \cup \tilde{U}^+_h({\Gamma_{\rm in}})$, where $(\vr,\vu)(t,x) = (\tilde\vr_B,\tilde\vu_B)$ on $\tilde{U}^+_h({\Gamma_{\rm in}})$, satisfies continuity equation in the sense of distributions on  $\Omega_h$ and 
$\vr \in L^2(0,T;L^2(\Omega_h))$ and $\vu \in L^2(0,T; W^{1,2}(\Omega_h))$. Then by classical DiPerna and Lions arguments with Friedrichs lemma provide that
\bFormula{rcebis}
\begin{split}
&\int_{\Omega_h} b(\vr)(\tau,\cdot)\varphi(\tau,\cdot) \,\dx - \int_{\Omega_h} b(\vr_0)(\cdot)\varphi(0,\cdot) \,\dx \\
&=
\int_0^\tau\int_{\Omega_h} \Big(b(\vr)\partial_t\varphi + b(\vr )\vu \cdot \Grad \varphi+(b(\vr)-b'(\vr)\vr)\Div\vu\Big) \dxdt 
\end{split}
\eF
holds for any $b\in C[0,\infty)$ with $b'\in C_c[0,\infty)$, $\tau\in [0,T]$, and $\varphi \in C_c^1([0,T]\times{\Omega_h})$.

In order to find boundary term $\int_0^\tau\int_{\Gamma_{\rm in}} b(\vr_B) \vu_B \cdot \vc{n} \varphi \ {\rm d}S_x{\rm d} t$ we write 
$$\int_{\Omega_h} b(\vr )\vu \cdot \Grad \varphi \,\dx = 
\int_{\Omega} b(\vr )\vu \cdot \Grad \varphi \,\dx + 
\int_{\tilde{U}^+_h({\Gamma_{\rm in}})} b(\vr )\vu \cdot \Grad \varphi \,\dx .
$$ 
By \eqref{ggg} and integration by parts
$$\int_{\tilde{U}^+_h({\Gamma_{\rm in}})} b(\vr )\vu \cdot \Grad \varphi \,\dx
= \int_{{\Gamma_{\rm in}}} b(\vr_B )\vu_B \cdot \vc{n} \varphi \,\dx
+
\int_{\tilde{U}^+_h({\Gamma_{\rm in}})} (\tilde\vr_B b(\tilde\vr_B ) - b(\tilde\vr_B)\Div \tilde\vu_B \,\dx.
$$
Inserting two above identities to \eqref{rcebis}, letting $h\to 0$, recalling regularity of $(\tilde\vr_B,\tilde\vu_B)$ we obtain desired conclusion of Lemma~\ref{Renorm}.

\smallskip

Furthermore, since \eqref{continuity} is satisfied also for $Z$ instead of $\vr$, due to \eqref{uweak} and by arguments of Lemma~\ref{Renorm}, $(Z,\vu)$ satisfies \eqref{P3}. That finishes the proof of Theorem~\ref{TM1!}. $\Box$

\subsection{Recovery of the system in terms of $(\vr,\vu,\vrs)$}\label{Sec:Recovery}

Our aim now is to prove that solution $(\vr,\vu, Z)$ can be identified with the solution $(\vr,\vu,\vrs)$  to the problem \eqref{main}. Namely, we need to show the existence of $\vr^*$ satisfying Definition~\ref{def_in_vrs}. 
To this end we will use combination of arguments from \cite[Section~3.1]{CHJN}, and from \cite[Section~4]{DeMiZa}.
First note that since $Z_0>0$ we have
$$\frac{\vr_0}{Z_0}= \vr^*_0.$$
Moreover, recall that we already know that $\vr$ and $Z$ satisfy renormalized continuity equations.

As in previous section, let us construct set 
 $\tilde{U}^+_h(\Gamma_{\rm in})$ and let us extend continuity equations for $(\vr,\vu)$ and $(Z,\vu)$ on 
 $\Omega_h:=\Omega\cup\Gamma_{\rm in} \cup \tilde{U}^+_h({\Gamma_{\rm in}})$. In particular extension $\tilde\vu_B$, $\tilde\vr_B$, and $\tilde{Z}_B$ of $\vu_B$, $\vr_B$, and ${Z}_B$ on $\tilde{U}^+_h(\Gamma_{\rm in})$ is, such that $\tilde\vu_B \in C^1(\tilde{U}^+_h(\overline{\Gamma_{\rm in})})$ and $\tilde\vr_B, \tilde{Z}_B \in W^{1,\infty}(\tilde{U}^+_h({\Gamma_{\rm in}}))$, and 
    \begin{equation}\label{ggg2}
    \begin{split}
   & \Div (\tilde\vr_B \tilde\vu_B ) = 0 \mbox{ in } \tilde{U}^+_h(\overline{\Gamma_{\rm in})}, \quad \tilde\vr_B|_{\Gamma_{\rm in}} = \vr_B, \  \tilde\vu_B|_{\Gamma_{\rm in}} = \vu_B,\\
   &
    \Div (\tilde{Z}_B \tilde\vu_B ) = 0 \mbox{ in } \tilde{U}^+_h(\overline{\Gamma_{\rm in})}, \quad \tilde{Z}_B|_{\Gamma_{\rm in}} = Z_B,
    \end{split}
    \end{equation}
Then extensions of $(\vr,\vu)$ and $(Z,\vu)$  to the set $\Omega_h:=\Omega\cup\Gamma_{\rm in} \cup \tilde{U}^+_h({\Gamma_{\rm in}})$, where $(\vr,\vu)(t,x) = (\tilde\vr_B,\tilde\vu_B)$ and 
$(Z,\vu)(t,x) = (\tilde{Z}_B,\tilde\vu_B)$ on $\tilde{U}^+_h({\Gamma_{\rm in}})$ 
satisfy continuity equations in the sense of distributions on  $\Omega_h$.

Applying convolution with a standard family of regularizing kernels we obtain the regularized functions $[\vr]_\omega$, $[Z]_\omega$  which satisfy 
\begin{equation}\label{reg_1omega}
    \partial_t [\vr]_\omega + \Div([\vr]_\omega\vu) = R^1_\omega \mbox{ a.e. in }(0,T)\times\Omega_{\omega,h}
\end{equation}
\begin{equation}\label{reg_2omega}
    \partial_t [Z]_\omega + \Div([Z]_\omega\vu) = R^2_\omega \mbox{ a.e. in }(0,T)\times\Omega_{\omega,h}
\end{equation}
where 
$$
\Omega_{\omega,h}: =
\{ x \in \Omega_h \, | \, {\rm dist}(x,\partial \Omega_h)>\omega \}.
$$
Due to Friedrichs commutator lemma, see e.g. \cite[Lemma~10.12]{FEINOV}, we find that 
$$
R^1_\omega \to 0 \mbox{ and } R^2_\omega \to 0 \mbox{ in } L^1_{\rm loc}((0,T)\times \Omega_h) \mbox{ as } \omega\to 0.
$$
Let us now multiply  \eqref{reg_1omega} by
$\frac{1}{[Z]_\omega+\lambda}$, and \eqref{reg_2omega}
by $- \frac{[\vr]_\omega + \lambda \vr^*_0}{([Z]_\omega+\lambda)^2}$,
with $\lambda>0$. Then after some algebraic manipulations we find that
\begin{equation*}
\begin{split}
\partial_t \left(
\frac{[\vr]_\omega + \lambda \vr^*_0}{[Z]_\omega+\lambda} \right)
+ 
\Div \left( \left( 
\frac{[\vr]_\omega + \lambda\vr^*_0}{[Z]_\omega+\lambda}
\right) \vu \right)
& - 
\left(
\frac{([\vr]_\omega + \lambda\vr^*_0)[Z]_\omega}{([Z]_\omega+\lambda)^2} 
+ \frac{\lambda \vr^*_0 }{[Z]_\omega+\lambda} 
\right) \Div \vu 
\\ & = R^1_\omega \frac{1}{[Z]_\omega+\lambda} 
- R^1_\omega \frac{[\vr]_\omega + \lambda \vr^*_0}{([Z]_\omega+\lambda)^2}
\quad \mbox{ a.e. in }(0,T)\times\Omega_{\omega,h}.
\end{split}
\end{equation*}
Testing above by $\varphi \in C_c^1([0,T]\times{\Omega_h})$, after integration by parts and after passing with $\omega \to 0$ we obtain that 
\bFormula{rcebiss}
\begin{split}
&\int_{\Omega_h} 
\left(
\frac{\vr + \lambda \vr^*_0}{Z+\lambda} \right)
(\tau,\cdot)\varphi(\tau,\cdot) \,\dx - \int_{\Omega_h} 
\left(
\frac{\vr_0 + \lambda \vr^*_0}{Z_0+\lambda} \right)
(\cdot)\varphi(0,\cdot) \,\dx \\
&=
\int_0^\tau\int_{\Omega_h} \Big(\left(
\frac{\vr + \lambda \vr^*_0}{Z+\lambda} \right)\partial_t\varphi 
+  \left( 
\frac{\vr + \lambda\vr^*_0}{Z+\lambda}
\right) \vu \cdot \Grad \varphi+
\left(
\frac{(\vr + \lambda\vr^*_0)Z}{(Z+\lambda)^2} 
+ \frac{\lambda \vr^*_0 }{Z+\lambda} 
\right)
\Div\vu\Big) \dxdt 
\end{split}
\eF
for any $\varphi \in C_c^1([0,T]\times{\Omega_h})$.

Next we distinguish two cases:\\
Case 1. For $Z=0$, due \eqref{fs-} we notice that  that $\vr = 0$ and therefore $\frac{\vr + \lambda\vr^*_0}{Z+\lambda} = \vr^*_0$ and 
$
\frac{(\vr + \lambda\vr^*_0)Z}{(Z+\lambda)^2} 
+ \frac{\lambda \vr^*_0 }{Z+\lambda} 
= \tilde\vr^*_0 $, then \eqref{rcebiss} becomes trivial.\\
Case 2. For $Z>0$, we find that  
$\frac{\vr + \lambda\vr^*_0}{Z+\lambda} \leq \max\{\vr^*_0, \frac{1}{c_*}\}.$
Since $\vr + \lambda\vr_0^* $ converges strongly  to $\vr$ as $\lambda \to 0$, as well as $Z + \lambda $ converges strongly  to $Z$ as $\lambda \to 0$, 
we can pass with $\lambda \to 0$ in \eqref{rcebiss} using Legesgue's Dominated convergence theorem,
to obtain that
\bFormula{rcebisss}
\begin{split}
&\int_{\Omega_h} \frac{\vr}{Z}(\tau,\cdot)\varphi(\tau,\cdot) \,\dx - \int_{\Omega_h} \frac{\vr_0}{Z_0}(\cdot)\varphi(0,\cdot) \,\dx \\
&=
\int_0^\tau\int_{\Omega_h} \Big(\frac{\vr}{Z}\partial_t\varphi + \frac{\vr}{Z}\vu \cdot \Grad \varphi+\frac{\vr}{Z}\Div\vu\Big) \dxdt 
\end{split}
\eF
holds for any $\varphi \in C_c^1([0,T]\times{\Omega_h})$. By the same steps as in  the proof of Lemma~\ref{Renorm} we find that after passing with $h\to 0$ we obtain 
\bFormula{rcebix}
\begin{split}
&\int_{\Omega} \frac{\vr}{Z}(\tau,\cdot)\varphi(\tau,\cdot) \,\dx - \int_{\Omega} \frac{\vr_0}{Z_0}(\cdot)\varphi(0,\cdot) \,\dx \\
&=
\int_0^\tau\int_{\Omega} \Big(\frac{\vr}{Z}\partial_t\varphi + \frac{\vr}{Z}\vu \cdot \Grad \varphi+\frac{\vr}{Z}\Div\vu\Big) \dxdt 
 - \int_0^\tau\int_{\Gamma_{\rm in}} \frac{\vr_B}{Z_B} \vu_B \cdot \vc{n} \varphi \ {\rm d}S_x{\rm d} t 
\end{split}
\eF
Obviously, $\vr^*$ defined as $\frac{\vr}{Z}$ satisfies $\vr^* \in \{ \min \{ \frac{1}{c^*} , \vr^*_0 \}, \max \{ \frac{1}{c^*} , \vr^*_0 \} \}$ a.e. in $(0,T)\times \Omega$, and thus $Z = \frac{\vr}{\vr^*}$ a.e. in $(0,T)\times \Omega$.

Finally we can conclude that Theorem~\ref{TM2} is proven. $\Box$

\section{Passage to the limit $\ep\to 0$} \label{Sec:lim}
The purpose of this section is to perform the limit $\ep\to0$ in the auxiliary system (\ref{main_transformed}-\ref{ia4}) to prove the Theorem \ref{TM3}. From now on, by $\{\vr_\ep, Z_\ep,\vu_\ep\}_{\ep>0}$ we denote the sequence of solutions obtained in the previous section.

\subsection{Convergence following from the uniform estimates}
The energy inequality established in the previous section, gives rise to the following estimates that are uniform with respect to $\ep$:
\begin{equation} \label{un_ep}
\begin{gathered}
\sup_{t\in[0,T]}\lr{\|\sqrt{\vr_\ep} \vu_\ep (t)\|_{L^2(\Omega)} 
+\left\|H(Z_\ep)(t) \right\|_{L^1(\Omega)} }\leq C,\\
\intT{{\|\vu_\ep\|_{W^{1,2}(\Omega,\R^3)}^2}} \leq C.
\end{gathered}
\end{equation}
From here and from the construction we also deduce that
\eq{\label{Zr_ep}
0\leq Z_\ep<1,\quad \mbox{a.e. in } Q_T,\quad 0\leq c_\star \vr_\ep\leq Z_\ep\leq c^\star\vr_\ep,}
in particular both sequences $Z_\ep$, $\vr_\ep$ are uniformly bounded in $L^p(Q_T)$ for any $p\leq \infty$. Therefore, up to the subsequence, we have
\eq{\label{conv_ep}
 &\vu_\ep \rightarrow \vu \qquad \text{weakly in }  L^2(0,T;W^{1,2}(\Omega, \R^3)),\\
&Z_\ep\rightarrow Z \qquad \text{in }\ C_{\rm weak}([0,T];L^{p}(\Omega)),\\
&\vr_\ep\rightarrow \vr \qquad \text{ in }\ C_{\rm weak}([0,T];L^{p}(\Omega)),
}
for any $p$ finite. In addition to that, the limiting $Z$ and $\vr$ will satisfy
\eq{\label{Zr_lim}
0\leq Z\leq1,\quad \mbox{a.e. in } Q_T,\quad 0\leq c_\star \vr\leq Z\leq c^\star\vr.}

To obtain the uniform $L^1$ bound for the pressure we can follow one of two strategies:

\noindent(i) for the zero flux, i.e. for $\int_{\partial \Omega} \vu_B \cdot \vc{n} \ {\rm d}S_x = 0$, we use our smallness assumption
\eq{
\intO{Z_0}+T\int_{\Gamma_in}Z_B|\vu_B\cdot\vc{n}|\ {\rm d}S_x<|\Omega|,
}
to see that it implies the condition \eqref{lambda} from the level of uniform in $\delta$ estimates, and we repeat the argument from the previous section.\\

\noindent(ii) for the positive flux, thanks to the energy estimate \eqref{ei!+} and the property \eqref{Lep} we already control the $L^1$ norm of $\pi_\ep$ uniformly on the subset ${\cal O}$. 
In order to control the pressure on the complementary part, we repeat the Bogovskii type estimate with the test function: 
\eq{\label{testB}
\bm{\psi}=\eta(t){\cal B}\lr{\xi},}
where
\eq{
\xi=\left\{\begin{array}{ll}
1& \mbox{in\ }\omega\setminus {\cal O}, \\
-\frac{|\Omega\setminus {\cal O}|}{|\Omega|} &\mbox{otherwise.}
\end{array}
\right.}
Note that since ${\cal O}$ is an open subset of $\Omega$, the test function $\bm{\psi} $is well defined and it belongs to $W^{1,p}_0(\Omega)$ for any $p<\infty$, so it is an admissible test function.

Similarly as in the previous limit passage we can in fact show that the local-in-time pressure bound holds true with $\eta \equiv 1$. However, we do not have anymore the equi-integrability of the pressure as on the previous level of approximation. Therefore, the convergence in the sense of measures is the most we can hope for, we have
\eq{\label{lim_pi_ep}
&\pi_{\ep}(Z_\ep) \rightarrow\pi  \quad\text{weakly\ in }\quad {\cal M}^+([0,T]\times K),\\
&Z_\ep \pi_{\ep}(Z_\ep) \rightarrow\pi_1  \quad\text{weakly\ in }\quad {\cal M}^+([0,T]\times K).
}
The limiting momentum equation therefore reads
 $$ \partial_t (\vr\vu) + \Div (\vr\vu \otimes \vu) + \Grad\pi -\Div\mathbb{S}(\Grad \vu) = \vc{0},$$
in the sense of distributions.

Already at this point we can identify the second limit in \eqref{lim_pi_ep} using the explicit form of the pressure \eqref{i5a}. We have
\eq{
Z_\ep\pi_\ep(Z_\ep)=\pi_\ep(Z_\ep)-\ep\frac{1}{(1-Z_\ep)^{\beta-1}},
}
thus letting $\ep\to0$ and observing that by \eqref{un_ep}, the last term converges to zero strongly, we obtain the relation
\eq{\label{step1}
\pi_1 =\pi  }
in the sense of distributions. 
The thing that remains to be shown in order to deduce the constraint $(1-Z)\pi=0$ is that we have $\pi_1=Z\pi$ in some sense (note that on the l.h.s. we have multiplication of measure by the $L^\infty$ function only).
The recovery of the constraint  requires stronger information about the convergence of $Z_\ep$, and additional information about the regularity of $\pi$ and $Z$.

\subsection{Strong convergence of $Z_\ep$}
We first need to show that we have a variant of effective viscous flux equality. It can be derived 
by testing the approximate momentum equation by the inverse divergence operator $\psi\eta\Grad\lap^{-1}[1_\Omega Z_\ep]$ and by testing the momentum equation by $\psi\eta\Grad\lap^{-1}[1_\Omega  Z]$, and by comparison of the limits, for any $ \psi \in C_c^\infty((0,T)) $ and $\eta \in C_c^\infty(\Omega)$.
Note that already at this stage we need to justify what the product $\pi Z$ means, i.e., whether 
$\eta\psi\Grad\lap^{-1}[1_\Omega Z]$ is regular enough to be used as a test function in the limiting momentum equation. 
To justify this step we first write weak formulation of the limiting momentum equation 
\eq{\label{pixi}
\langle\pi,\Div\xxi\rangle_{({\cal M}(Q_T),C(Q_T))}
&=\intTO{\mathbb{S}(\Grad \vu):\Grad\xxi} 
-\intTO{\vr\vu\cdot\pt\xxi}\\
&\quad -\intTO{\vr\vu\otimes\vu:\Grad\xxi} -\intTO{\vr(\vc{w}-\vu)\cdot \xxi}}
that is satisfied for all $\xxi\in C^1_c(Q_T)$. From now on we will treat this formulae as a definition of $\intTO{\pi\Div\xxi}$. Let us check that the r.h.s. of \eqref{pixi} makes sense for $\xxi$ from much wider class, it is enough that
\eqh{
\Grad \xxi\in L^2(0,T; L^2(\Omega;R^{d\times d})),\quad\Grad\xxi\in L^{5/2}(0,T; L^{5/2}(\Omega;R^{d\times d})),\quad \pt\xi\in L^1(0,T; L^2(\Omega;R^d)).
}
The second property is a consequence of simple interpolation property
$$\|\vr|\vu|^2\|_{L^{5/3}(Q_T)}\leq
\|\vr|\vu|^2\|_{L^\infty(0,T; L^1(\Omega))}^{2/5}\|\vr|\vu|^2\|_{L^1(0,T; L^3(\Omega))}^{3/5}.
$$
So, if we take
\eq{\label{class_xi}
\xi\in W^{1,5/2}_0(Q_T;R^d),}
 the r.h.s. of \eqref{pixi} will be well defined. Note that our $C^1_c(Q_T)$ functions are dense in $W^{1,5/2}_0(Q_T)$.

Let us check that $\xxi:=\psi\phi\Grad\lap^{-1}[1_\Omega Z]$ has these properties.
First, note that the Riesz operator 
$$\A=\Grad\lap^{-1}: \ L^p(R^d)\to D^{1,p}(R^d;R^d)$$ (homogeneous Sobolev space)
is a continuous linear operator and we have that
$$\|\Grad\A[v]\|_{L^p(R^d;R^d)}\leq C(p)\|v\|_{L^p(R^d)},$$
for any $1<p<\infty$.
Therefore
\eqh{\|\Grad\xxi\|_{L^\infty(0,T;L^p(R^d;R^{d\times d}))}&=
\|\Grad\lr{\psi\phi\Grad\lap^{-1}[1_\Omega Z]}\|_{L^\infty(0,T;L^p(R^d;R^d))}\\
&\leq C(p,\psi,\phi) (1+\|1_\Omega Z\|_{L^\infty(0,T;L^p(R^d))})\leq C,}
for any $p<\infty$.
Next, using the continuity equation we get that
\eqh{
\pt\xxi=\pt\lr{\psi\phi\Grad\lap^{-1}[1_\Omega  Z]}&=\phi\pt\psi\Grad\lap^{-1}[1_\Omega  Z]+\phi\psi\Grad\lap^{-1}[1_\Omega  \pt Z]\\
&=\phi\pt\psi\Grad\lap^{-1}[1_\Omega Z]-\phi\psi\Grad\lap^{-1}[ \Div(1_\Omega Z\vu)].
}
Using the properties of the double Riesz transform we obtain 
\eq{
\|\pt\xxi\|_{L^p(0,T;L^p(\Omega))}\leq C(p,\psi,\phi)(1+ \|Z\vu\|_{L^p(0,T;L^p(\Omega))})\leq C,
}
for some $p>5/2$. Whence, we have shown that $\xxi=\psi\phi\Grad\lap^{-1}[1_\Omega Z]$ is indeed in the class \eqref{class_xi}.

Now, note that 
$$\Div\xxi=\psi\Grad\phi\cdot\Grad\lap^{-1}[1_\Omega Z]+\psi\phi Z,$$
so we can define
\eq{
&\langle\pi,\phi\psi Z\rangle_{({\cal M}(Q_T),C(Q_T))}:=\langle\pi,\Div\xxi\rangle_{({\cal M}(Q_T),C(Q_T))}-
\langle\pi,\psi\Grad\phi\cdot\Grad\lap^{-1}[1_\Omega Z]\rangle_{({\cal M}(Q_T),C(Q_T))}.
}
This means that
$\langle\pi,\phi\psi Z\rangle_{({\cal M}(Q_T),C(Q_T))}$ is well defined iff 
$\langle\pi,\phi\psi \Grad\lap^{-1}[1_\Omega Z]\rangle_{({\cal M}(Q_T),C(Q_T))}$
 is well defined. For that we need $ \Grad\lap^{-1}[1_\Omega Z]$ to be at least $C(Q_T)$, but this is true, as $Z\in C_{\rm weak}(0,T; L^p(\Omega))$ for sufficiently large $p$. In particular, for $p>d$, using the Morrey inequality, we get $ \Grad\lap^{-1}[1_\Omega Z]\in C([0,T]\times \Ov{\Omega};R^d)$.

Taking the above into account, we have
 
\eq{ \label{evf_ep}
&\lim_{\ep \to 0^+}  \intTO{\psi  \phi \big( \pi_{\ep}(Z_\ep)-(\lambda+2\mu) \Div \vu_\ep\big) Z_\ep}\\
&=  \langle\pi,\phi\psi Z\rangle_{({\cal M}(Q_T),C(Q_T))}-(\lambda+2\mu) \intTO{\psi \phi    \Div \vu Z }.
}
From \eqref{evf_ep} it follows that
\eq{\label{lim}
&(\lambda+2\mu)\intTO{\psi\phi\lr{\Ov{Z\Div\vu}-Z\Div\vu}}\\
&=\langle\pi_1,\phi\psi \rangle_{({\cal M}(Q_T),C(Q_T)}
-\langle\pi,\phi\psi Z\rangle_{({\cal M}(Q_T),C(Q_T))}\\
&=\langle\pi,\phi\psi (1-Z)\rangle_{({\cal M}(Q_T),C(Q_T))}\geq 0
}
where we have used subsequently \eqref{step1}, and the limit of \eqref{Zr_ep}. Since both pairs $(Z_\ep,\vu_\ep)$ and $(Z,\vu)$ satisfy the renormalized continuity equation, we can use the renormalization in the form $b(z)=z\log z$ to justify that 
$$Z_\ep\to Z\quad \text{strongly in } L^p((0,T)\times \Omega),\quad \forall p<\infty.$$
The proof is identical as for the limit passage $\delta\to 0$.

Note, however, that similarly as in the previous section this property is not transferred to the sequence $\vr_\ep$, for which we only have \eqref{conv_ep}. 
Nevertheless, using this information and formula \eqref{lim} we can justify that
$$\langle\pi,\phi\psi (1-Z)\rangle_{({\cal M}(Q_T),C(Q_T))}=0.$$

Note that at this stage we can repeat the arguments from Section \ref{Sec:Recovery} in order to come back to the solution in terms of $\vr,\vu,\vr^*$. Indeed, note that the proof is based only on the properties of the renormalized continuity equations for $\vr$ and $Z$, on the boudedness of $Z$ and $\vr$, and on the regularity of $\vu$ which are the same as on the previous level of approximation. 

Having this, justification of condition \eqref{cond:u} amounts to repetition of proof of \cite[Lemma 4]{PeZa}, see also \cite[Lemma 2.1]{LM99}. The proof of Theorem \ref{TM3} is thus complete. $\Box$



	
\def\cprime{$'$} \def\ocirc#1{\ifmmode\setbox0=\hbox{$#1$}\dimen0=\ht0
  \advance\dimen0 by1pt\rlap{\hbox to\wd0{\hss\raise\dimen0
  \hbox{\hskip.2em$\scriptscriptstyle\circ$}\hss}}#1\else {\accent"17 #1}\fi}

\end{document}